\documentclass[onecolumn,11pt]{article}

\usepackage[a4paper,margin=1in]{geometry}

\usepackage[affil-it]{authblk}

\usepackage[english]{babel}

\usepackage[onehalfspacing]{setspace}

\usepackage{amsmath,amssymb,amsthm,amsfonts,mathtools,mathrsfs}
\usepackage{graphics}
\usepackage{subcaption}
\usepackage{bm}
\usepackage{refcount}
\usepackage{booktabs}
\usepackage{enumerate}
\usepackage{adjustbox}
\usepackage{pifont}

\usepackage[short,nocomma]{optidef}

\usepackage[utf8]{inputenc} 

\usepackage{amsthm}
\newtheorem{theorem}{Theorem}

\newtheorem{lemma}[theorem]{Lemma}
\newtheorem{corollary}[theorem]{Corollary}
\newtheorem{proposition}[theorem]{Proposition}

\newtheorem{definition}[theorem]{Definition}

\usepackage[authoryear,round]{natbib}
\usepackage{algorithm}
\usepackage[noend]{algpseudocode}

\DeclareMathOperator*{\argmin}{argmin}

\newcommand{\timetable}{\theta}
\newcommand{\timetableBold}{{\bm{\theta}}}

\newcommand{\routes}{\mathcal{R}}

\newcommand{\routingBold}{\bm{p}}
\newcommand{\xmark}{\ding{55}}
\newcommand{\cmark}{\ding{51}}

\usepackage[colorinlistoftodos]{todonotes}

\usepackage{xcolor,soul}
\definecolor{dr}{rgb}{0.75,0.00,0.00}
\definecolor{lr}{rgb}{1.00,0.75,0.75}
\sethlcolor{lr}

\usepackage{graphicx}
\usepackage{tikz}
\usepackage{pgfplots}
\usetikzlibrary{pgfplots.statistics, pgfplots.colorbrewer} 
\pgfplotsset{compat=1.18} 
\usetikzlibrary{shapes,arrows,positioning}
\tikzstyle{vertex} = [draw,rounded rectangle,thick]

\usepackage[
	colorlinks = true,
	linkcolor = blue,
	urlcolor  = blue,
	citecolor = black,
	anchorcolor = blue
]{hyperref}

\DeclareRobustCommand{\VAN}[3]{#2} 

\makeatletter
\def\@maketitle{%
	\newpage
	\null
	\vskip 2em%
	\begin{center}%
		\let \footnote \thanks
		{\Large\bfseries \@title \par}%
		\vskip 1.5em%
		{\normalsize
			\lineskip .5em%
			\begin{tabular}[t]{c}%
				\@author
			\end{tabular}\par}%
		\vskip 1em%
		{\normalsize \@date}%
	\end{center}%
	\par
	\vskip 1.5em}
\makeatother

\title{A Unified Approach to Evaluation and Routing in Public Transport Systems}

\author[1]{Rolf Nelson van Lieshout\footnote{Corresponding author: r.n.v.lieshout@tue.nl}}
\author[2]{Kevin Dalmeijer}

\affil[1]{Department of Operations, Planning, Accounting, and Control,\protect\\School of Industrial Engineering, Eindhoven University of Technology\vspace{\baselineskip}}
\affil[2]{H. Milton Stewart School of Industrial and Systems Engineering,\protect\\Georgia Institute of Technology}

\date{\today}

\graphicspath{{./figures/}}

\begin{document}

\maketitle

\thispagestyle{empty}

\vfill

\begin{abstract}
\noindent \textbf{\textit{Problem definition:}} Both evaluating the service quality of a public transport system and understanding how passengers choose between modes or routes is imperative for public transport operators, providers of competing mobility services and policy makers. However, the literature does not offer consensus on how either of these tasks should be performed, which can lead to inconsistent or counter-intuitive results.  \textbf{\textit{Methodology/results:}} This paper provides a formal treatment on how fundamental elements of public transport systems (route sets, timetables and line plans) can be evaluated consistently, and how travelers distribute over routes. Our main insight is that evaluation and routing are two sides of the same coin: by solving an appropriate optimization model one obtains both the quality of the route set, timetable or line plan (the optimal objective value), and the distribution of the travelers over the routes (the optimal solution itself). The practical relevance of the new framework is demonstrated with several applications that are validated with real data from the Dutch and Swiss railway networks. \textbf{\textit{Managerial implications:}} The measures and route choice models developed in this paper enable planners to create better line plans and to effectively analyze timetables for inefficiencies. The framework also reveals: (i) the importance of using the right model for the right stage of planning, (ii) that it is not always necessary for public transport planners to accurately model travel behavior, especially for high-level planning, and (iii) that combining models in an inconsistent way can have significant negative consequences that are avoided with the new framework.\\

\end{abstract}

\clearpage

\section{Introduction}
\label{sec:intro}
In passenger mobility research, it is common to define some measure for the service quality of a public transport system (e.g., average travel time) and to assume some choice model that distributes travelers over routes and/or modes (e.g., the multinomial logit model). Despite their ubiquity, there is currently no consensus on \textit{how} the service quality of public transport should be measured and \textit{what} choice model is appropriate for what context. Instead, the literature offers a variety of approaches, which makes it difficult to compare and validate obtained results. Furthermore, it is generally accepted that public transport systems can improve service quality by adding routes, making routes faster or increasing the frequency at which routes are operated (at least, in the absence of congestion effects). However, this is not always reflected in existing approaches, leading to counter-intuitive and inconsistent results. 

\begin{figure}[bh]
    \centering
	    \begin{subfigure}[t]{0.165\textwidth}
			\centering
	       	\includegraphics[width=0.8\textwidth, trim= 0 -2cm 0 0, clip]{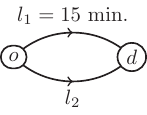}
				\caption{Routes.}
        \label{fig:paradox:a}
		\end{subfigure}
		\begin{subfigure}[t]{0.41\textwidth}
			\centering
	       	\includegraphics[width=\textwidth]{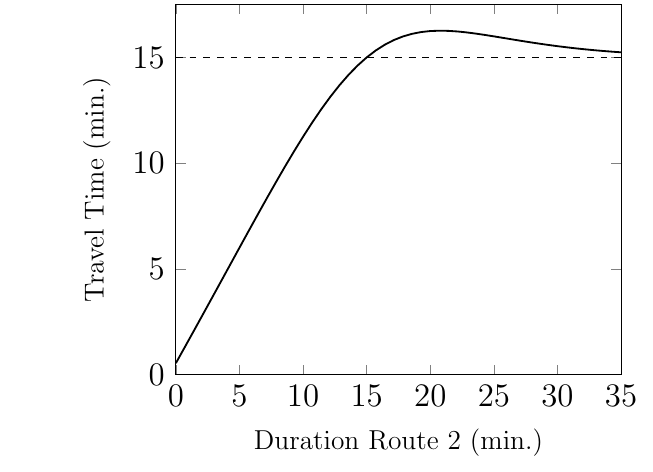}
				\caption{Travel Time as function of $l_2$.}
			\label{fig:paradox:b}
		\end{subfigure}
		\begin{subfigure}[t]{0.41\textwidth}
			\centering
	       	\includegraphics[width=\textwidth]{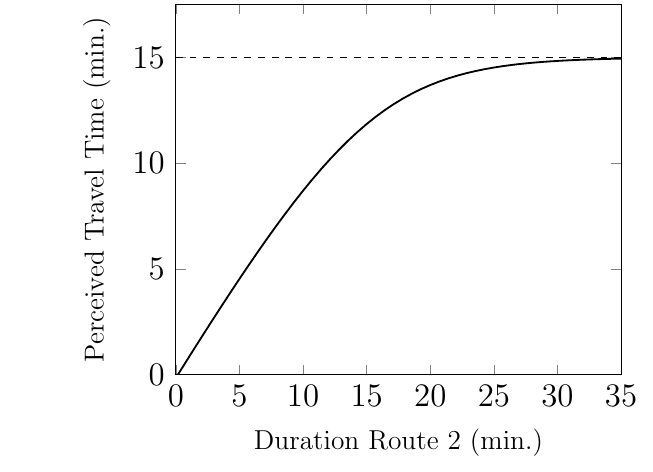}
				\caption{Perceived Travel Time as function of $l_2$.}
			\label{fig:paradox:c}
		\end{subfigure}
    \caption{Illustration of two route set measures, where travelers are distributed over the routes using the logit model with $\beta=0.22$. The duration of route $i$ is denoted by $l_i$.}
    \label{fig:paradox}
\end{figure}

To illustrate how seemingly sensible route choice models and service quality measures can lead to unexpected outcomes and suboptimal decisions, suppose that travelers choose routes according to the logit model and that average travel time is used as a measure of service quality. For the system depicted in Figure~\ref{fig:paradox:a}, where travelers can choose between route 1 with a duration of 15 minutes and route 2 with duration of $l_2$ minutes, Figure~\ref{fig:paradox:b} shows the average travel time as a function of $l_2$. After some point, making route 2 slower actually \emph{improves} the measure, because it increases the likelihood that travelers switch to route 1, which is better in terms of travel time. Public transport planners using such a measure are insufficiently incentivized to speed up or add routes. In fact, they may even decide to slow down or remove routes, as it could occur that their measure suggests that doing so improves service quality. Note that this phenomenon is fundamentally different than Braess' paradox, because in our setting travel times are fixed, whereas in Braess' paradox travel times depend on the number of travelers that choose a route \citep{Braess1968-UberEinParadoxon}. 

\begin{figure}[t]
		\centering
		\begin{subfigure}[t]{0.19\textwidth}
			\centering
			\fbox{\includegraphics[height=2.3cm]{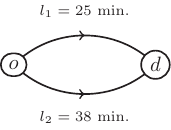}}
			\caption{Route Set.}
            \label{subfig:rs}
		\end{subfigure}
		\hspace{0.01\textwidth}
		\centering
		\begin{subfigure}[t]{0.372\textwidth}
			\centering
			\fbox{\includegraphics[height=2.3cm]{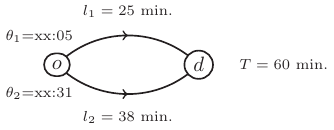}}
			\caption{Timetable.}
            \label{subfig:tt}
		\end{subfigure}
		\hspace{0.01\textwidth}
		\begin{subfigure}[t]{0.372\textwidth}
			\centering			\fbox{\includegraphics[height=2.3cm]{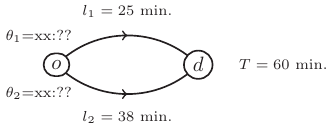}}
			\caption{Line Plan.}
            \label{subfig:lp}
		\end{subfigure}
		\caption{Differences between a route set, timetable and line plan. $l_i$ and $
		\theta_i$ represent the route duration and departure time of route $i$, respectively. $T$ is the cycle time or \textit{period}.}
		\label{fig:concepts}
\end{figure}

In this paper, we present a formal treatment of route choice models and service quality measures for three supply models of public transport: route sets, (periodic) timetables and line plans, concepts illustrated in Figure~\ref{fig:concepts}.
The theory is built for a single origin-destination (OD) pair, leading to a framework that can be applied at the network level.
We analyze the two predominant route choice models in the literature: \textit{shortest path} routing, where all travelers choose the shortest route, and \textit{logit} routing, where travelers distribute over the routes according to the logit model. We then define desirable properties of measures, show which measures fail to meet these properties and develop measures that do. To ensure that the derived measures are consistent and interpretable, our approach is hierarchical: the line plan measures build upon the timetable measures, which again build upon the route set measures. Ultimately, this results in a ready-to-use framework for routing and evaluation in public transport, and several applications are presented to demonstrate its practical relevance.

\paragraph{Theoretical Framework}

The framework will be developed starting from the concept of a \emph{route set}: a set of routes with given durations, as visualized by Figure~\ref{subfig:rs}.
A measure for route set quality can be constructed by pairing a routing model, e.g.\ logit, with an evaluation function, e.g.\ travel time.
Two desirable properties will be defined: \textit{monotonicity} and \textit{consistency}. A measure is monotonic if increasing route durations or removing a route cannot improve the measure.
Figure~\ref{fig:paradox:b} illustrates that the combination of logit routing and travel time evaluation induces a measure that fails to be monotonic.
We say a measure is consistent if the routing model minimizes the evaluation function.
Again, the combination of logit routing and travel time evaluation fails to satisfy this property, since travel time is minimized by shortest path routing, not by logit routing.

The observation that logit routing does not minimize travel time raises the question if there exists an alternative evaluation function that is minimized by logit. We answer this question affirmatively, by showing that logit routing optimizes the evaluation function that we refer to as \textit{perceived travel time}, the sum of travel time and an additional entropy term. This results in a measure that is both consistent and monotonic. The monotonicity can also be observed in Figure~\ref{fig:paradox:c}, which depicts perceived travel time under logit routing. Moreover, we provide a motivation for the terminology ``perceived travel time'' through the random utility model (RUM) interpretation of logit. In addition to evaluating perceived travel time under logit routing, combining regular travel time evaluation with shortest path routing also induces a consistent and monotonic measure. 

To also support practitioners in other stages of the public transport planning process, this paper subsequently extends the analysis from route sets to periodic timetables. A \emph{periodic timetable} is a route set that is operated periodically with a given cycle time or period and with given departure times (see Figure~\ref{subfig:tt}). 
The quality of a route now not only depends on its duration, but also on the waiting (or adaptation) time a traveler experiences: how much the departure time of a route deviates from their \textit{preferred} departure time \citep{friedrich2021departure}. If the period of the timetable is reasonably small, it is natural to assume that travel demand is uniformly distributed over the period \citep{gentile2016modelling,Hartleb2022,kaspi2013service,polinder2021timetabling}. Under this assumption, we show how to efficiently compute the average travel time and the average perceived travel time (both including waiting time), resulting in two consistent and monotonic measures for the quality of a periodic timetable.  

Finally, we consider \emph{line plans}, which are route sets that will be performed periodically, but with unknown departure times (see Figure~\ref{subfig:lp}). Line plans are commonly used in strategic planning, when the timetable is not yet known. There, one needs measures for service quality and route choice models that are independent of the timetable that will be operated. To deal with this additional source of uncertainty, we construct line plan measures by optimizing the respective timetable measures over all possible timetables. We additionally develop algorithms that efficiently solve this optimization problem both for shortest path and for logit routing.

The developed framework reveals a deep parallel between evaluation and routing: by solving an appropriate optimization model one obtains both the quality of the route set, timetable or line plan (the optimal objective value), and the distribution of the travelers over the routes (the optimal solution itself). Hence, every evaluation function \textit{implies} some route choice model, and vice versa. The importance of this connection is evident for route sets from the example in Figure~\ref{fig:paradox}, but this behavior extends all the way to line plans: it will be shown that if the evaluation function and routing model do not line up, it can lead to missed opportunities for improving services and wrong decisions that can harm the experience of public transport users. 

\paragraph{Practical Applications and Managerial Insights}

To demonstrate the practical relevance of the new framework, we present several applications that are validated with real data from the Dutch and Swiss railway networks.
We first focus on line planning, and we discuss how the new line plan measures can be integrated into existing methods to improve line planning practices.
To empirically validate the proposed line plan measures, we compare them to the timetable measures of the \emph{actual} timetable that was implemented in practice.
Remarkably, it is found that the line plan measures are 99\% accurate in predicting the measure of the real system, effectively anticipating the impact of a timetable that was not yet known at the line planning stage.
This suggests that the new line plan measures have major practical value in capturing the impact of downstream decisions at the line planning stage.
A similar analysis shows that 91\% of passenger flows can be predicted correctly at the line planning stage, opening the door for high-level capacity planning and congestion analysis without the need for timetables or simplifying assumptions.

Second, we show how the new framework can be used to automatically identify potential inefficiencies in a given timetable.
This is especially important when timetabling is partially manual and involves thousands of OD pairs, as may be the case in practice \citep{HuismanMaroti2024-OperationsResearchNetherlands}.
We apply this method to the Dutch network to identify an inefficient connection between Rotterdam and Bijlmer, and we show how it could potentially be improved.
It is also discussed how this kind of analysis can be used to study the fairness of the timetable between different regions or connections.

Third, we use the new framework to obtain managerial insights about 1) the importance of using the right model at the right time, 2) when shortest path routing is all you need, and 3) the consequences of using inconsistent models.
By zooming in on the connection between Tilburg and Eindhoven in the Netherlands, we demonstrate that route choice and service quality can differ substantially between route
sets, timetables, and line plans, highlighting that it is crucial to select the right model at the right stage of planning.
We also use real data to show that, while logit is the model of choice to accurately capture traveler behavior, the benefit of logit over shortest path diminishes as we move from route sets to timetables to line plans.
This observation has major managerial implications, as it implies that \emph{it is not always necessary for public transport planners to accurately model travel behavior}. Especially for long-term or strategic planning such as line planning,
the shortest path model may be all you need to make good decisions.
Finally, we use an example to demonstrate that combining the models in an inconsistent way can have significant negative consequences, and that being consistent may be more important than being correct.

\paragraph{Contributions and Structure of the Paper}

In summary, the contribution of this paper is twofold. Firstly, it develops a framework for traveler route choice and service quality evaluation for route sets, timetables and line plans. Compared to existing work, the framework introduces a new evaluation function that is consistent with the logit model, and routing models and measures tailored to the specific nature of line plans. Secondly, it demonstrates the practical relevance of the new framework with several applications that are validated with real data.

The remainder of this paper is structured as follows. Section~\ref{sec:lit} describes related literature. Sections~\ref{sec:rs}, \ref{sec:tt} and \ref{sec:lp} develop the framework for route sets, timetables and line plans, respectively. Section~\ref{sec:ill} presents applications and managerial insights. Section~\ref{sec:con} concludes the paper and describes directions for future research.

\section{Literature Review}
\label{sec:lit}

Existing measures of public transport service quality can broadly be divided into two categories. On one hand, there are \textit{travel time-based} approaches, which evaluate quality using (generalized) travel time as a proxy for user satisfaction. On the other hand, \textit{utility-based} approaches, rooted in discrete choice theory, quantify the total utility experienced by users. Although this distinction is not always made explicit, it mirrors the divide between commonly used route choice models: shortest path models and variants of the multinomial logit model.

\subsection{Travel Time-Based Evaluation and Shortest Path Routing}

Travel time-based measures are predominant in the Operations Research literature due to their linearity and compatibility with mathematical optimization models. Many recent contributions in line planning—where decisions are made about direct lines and frequencies—and timetabling—where departure and arrival times are determined—minimize travel time. Examples include \citet{bertsimas2021data} and \citet{YAO2024102989} for line planning, and \citet{borndorfer2017} and \citet{schiewe2020periodic} for timetabling.

While such models capture the tendency of travelers to prefer faster routes, they fail to account for the fact that having \textit{more} available routes (i.e., higher frequencies) also improves service quality. This can be incorporated by expanding the origin-destination matrix with a temporal component to model when people want to travel, and include waiting or adaptation time in the generalized travel time \citep{kaspi2013service,polinder2021timetabling,RobenekEtAl2016-PassengerCentricTrain}. In periodic timetabling models, uniform demand over the period is typically assumed, which naturally rewards more departures, or more uniformly spaced departures. For instance, a timetable with departures at xx:05 and xx:35 yields lower expected waiting times than one with departures at xx:05 and xx:20.

In line planning, incorporating waiting time is more challenging because it depends on the eventual timetable, which is not yet fixed at the planning stage. Again assuming that demand is uniformly distributed over the period, \citet{GKIOTSALITIS2022103492} estimate expected waiting times under the assumption of independent, random departure times. Conversely, \citet{Hartleb2022} assume equidistant departures and a uniform distribution of passengers over all services, resulting in expected waiting times of $T/(2n)$ if $n$ services are available per period.

However, these assumptions decouple the line planning and timetabling stages and disregard that planners will eventually optimize the timetable. Our framework addresses this gap by proposing line planning measures and route choice models that explicitly incorporate the fact that departure times will be optimized in a subsequent timetabling phase. In doing so, we allow waiting time to be endogenously minimized, and we develop efficient algorithms for computing the resulting service quality measures and routings.

Another limitation of these travel time-based approaches is that they assume that all travelers take the single best route, which fails to reflect observed traveler behavior. Empirical studies show that travelers often split across routes with similar (generalized) travel times \citep{van2014deduction}. This motivates the utility-based methods discussed next.

\subsection{Utility-Based Evaluation and Logit Routing}

Utility-based approaches, commonly used in transportation economics, model service quality from the perspective of individual traveler choices. These methods are not constrained by linear optimization requirements and are built on the foundations of discrete choice theory \citep{ben1985discrete,mcfadden1973conditional}. Here, service quality is often referred to as \textit{accessibility}, \textit{user benefits}, \textit{consumer surplus}, or \textit{welfare}.

The dominant routing model in this domain is the logit model, which naturally induces the \textit{logsum} as a measure of service quality. The logsum reflects the expected utility across all available options \citep{williams1977formation}, but its absolute value is not directly interpretable. As such, evaluations typically focus on changes in the logsum relative to a baseline scenario \citep{sweet1997aggregate}, which can be converted into monetary or temporal units for use in cost-benefit analysis \citep{cats2022beyond,de2007logsum,standen2019value}.

However, most utility-based models do not distinguish between route sets, timetables, and line plans. They treat services as static sets of alternatives and fail to account for the importance of timing and frequency in shaping service quality measures and routing decisions. Our proposed approach directly addresses this limitation by tailoring the evaluation and routing models to each of the three public transport supply models.

\subsection{Logit Routing in Optimization Models}

Recent studies have begun to incorporate advanced route choice models into optimization frameworks. For example, \citet{banerjee2024plan} jointly optimize pricing and network design under logit-based route choice, effectively maximizing the total expected utility (i.e., logsum). Other works, such as \citet{bertsimas2020joint} and \citet{HartlebSchmidt2022-RailwayTimetablingIntegrated}, apply logit routing but retain travel time as the objective function. This introduces a subtle but important inconsistency: if travelers do not take the shortest route, minimizing travel time is no longer aligned with their actual behavior, leading to the counterintuitive phenomenon described in Section~\ref{sec:intro}.

In this paper, we formalize this inconsistency and show that such models can lead to measures that violate desirable properties—namely, consistency and monotonicity. As an alternative, we propose an objective function—\textit{perceived travel time}—that is both consistent with logit routing and monotonic, and easy to interpret. Moreover, our framework generalizes across route sets, timetables, and line plans, and enables efficient integration into planning models.

\section{Framework for Routing and Evaluation}
This section introduces a novel and ready-to-use framework for routing and evaluation in public transport.
New routing models and evaluation measures are built up hierarchically from route sets to periodic timetables and all the way up to line plans.
Along the way, we define desirable properties, show which measures fail to meet these properties, and develop measures that do.
We introduce optimization models to efficiently evaluate the new measures, and show that these optimization models can be interpreted as routing models themselves, proving that evaluation and routing are two sides of the same coin.

\subsection{Route Sets}
\label{sec:rs}

We start building the new framework from the most basic object: a set of routes $\routes = \{1,2,\hdots,n\}$ that travelers can choose from.
Each route $i \in \routes$ has a duration $l_i \in \mathbb{R}$, which may include more general quality aspects (all converted to time) such as the number of transfers, transfer time, and waiting time.
Travelers are assumed to make probabilistic decisions, and assign probabilities to the routes according to a \emph{route choice model}.
When the probabilities are set, the expected quality of the route set can be evaluated with an \emph{evaluation function}.
We refer to the combination of a route choice model and an evaluation function as a \emph{measure} of route set quality.

\paragraph*{Route Choice Models}
This paper considers the two predominant route choice models in the literature: shortest path routing and logit routing.
Let a route choice model, or a \emph{routing} for short, be a function that assigns probabilities to the routes.
That is, a routing is a function $\routingBold$ that maps the routes $\routes$ to $\routingBold(\routes) \in \mathcal{P}$, where $\mathcal{P}$ is the set of probability vectors.
The two routing models are defined as follows:

{\singlespacing
\begin{itemize}
	\item \emph{Shortest Path Routing:}
	\begin{equation*}
		p_i^{sp} = 1 \textrm{ if } i = \argmin_{j\in \routes} \left\{l_j\right\} \textrm{ and } p_i^{sp} = 0 \textrm{ else, for all } i \in \routes.
	\end{equation*}
	\item \emph{Logit Routing:}
	\begin{equation*}
	    p_i^{logit} = \frac{e^{-\beta l_i}}{\sum_{j\in \routes} e^{-\beta l_j}}, \textrm{ for all $i \in \routes$, with parameter $\beta > 0$.}
	\end{equation*}
\end{itemize}
}

\noindent Shortest path routing assigns probability one to the shortest path and probability zero otherwise, where ties are broken arbitrarily (e.g., by lowest index).
Logit routing assigns probabilities accordingly to the multinomial logit model \citep{mcfadden1973conditional}.
The parameter $\beta$ controls the sensitivity of travelers to the route durations.
For smaller $\beta$, travelers distribute more evenly, while for larger $\beta$ they are more likely to choose the shortest route.
In fact, it can be seen that shortest path routing is the limiting case of logit routing for $\beta \rightarrow \infty$.
Figure~\ref{fig:routingModels} provides an example for $\routes=\{1,2\}$ and $l_1=15$ minutes.
For shortest path routing, travelers choose route 2 if the duration is shorter than 15 minutes, and route 1 otherwise.
For logit routing, the probability to choose route 2 is a smooth function of $l_2$, the steepness of which depends on $\beta$. 

\begin{figure}[t]
 \centering
 \def\scaleRouting{0.6}
 \def\t2{15}
\begin{subfigure}[b]{0.4\textwidth}
        \centering
        \begin{tikzpicture}[scale=\scaleRouting]
        \def\bSP{28.9}
        
        \begin{axis}[cycle list name=black white,xmin=0,xmax=30,samples=500,ylabel={Probability Route 2},xlabel={Travel Time Route 2 (min.)}]
            \addplot+[thick,mark=none,domain = 0:15] {1};
             \addplot+[thick,mark=none,domain = 15:30] {0};
             \draw[fill=white] (15,1) circle (2pt);
            \draw[fill=black] (15,0) circle (2pt);

        \end{axis}
        \end{tikzpicture}
\caption{Shortest Path Routing.}
\label{fig:probSP}
\end{subfigure}
        \begin{subfigure}[b]{0.4\textwidth}
        \centering
        \begin{tikzpicture}[scale=\scaleRouting]
        \def\bOne{0.05}
        \def\bTwo{0.1}
        \def\bThree{0.2}
        \def\bFour{0.4}
        \def\bFive{0.5}
        \def\bSix{1.0}
        
        \begin{axis}[cycle list name=black white,xmin=0,xmax=30,samples=500,ylabel={Probability Route 2},xlabel={Travel Time Route 2 (min.)}]
                \addplot+[thick,mark=none,domain = 0:30] {exp(-\bSix*x)/( exp(-\bSix*\t2)+exp(-\bSix*x) )};
                \addplot+[thick,white!20!black,mark=none,domain = 0:30] {exp(-\bFive*x)/( exp(-\bFive*\t2)+exp(-\bFive*x) )};
                \addplot+[thick,white!40!black,mark=none,domain = 0:30] {exp(-\bThree*x)/( exp(-\bThree*\t2)+exp(-\bThree*x) )};
                \addplot+[thick,white!60!black,mark=none,domain = 0:30] {exp(-\bTwo*x)/( exp(-\bTwo*\t2)+exp(-\bTwo*x) )};
                \addplot+[thick,white!80!black,mark=none,domain = 0:30] {exp(-\bOne*x)/( exp(-\bOne*\t2)+exp(-\bOne*x) )};

            \legend{$\beta=\bSix$,$\beta=\bFive$,$\beta=\bThree$,$\beta=\bTwo$,$\textcolor{white}{s}\beta=\bOne$} 
        \end{axis}
        \end{tikzpicture}
\caption{Logit Routing.}
\label{fig:probL}
\end{subfigure}
\caption{Routing models for route set $\routes=\{1,2\}$ with $l_1=15$ minutes and varying $l_2$.}
\label{fig:routingModels}
\end{figure}

\paragraph*{Evaluation Functions}
After the route probabilities $p \in \mathcal{P}$ have been determined, the quality of the route set is evaluated with an evaluation function $\mathcal{E}(\routes, p) \in \mathbb{R}$ that assigns lower scores to better outcomes.
We consider two different evaluation functions: the expected \emph{travel time} and a new \emph{perceived travel time}.

\begin{itemize}
	\item \emph{Travel Time:} $$\mathcal{E}_{tt}(\routes, p) = \sum_{i \in \routes} l_i p_i.$$
	\item \emph{Perceived Travel Time:} $$\mathcal{E}_{ptt}(\routes, p) = \sum_{i \in \routes} l_i p_i + \frac{1}{\beta}\sum_{i\in \routes} p_i \log(p_i), \textrm{ with parameter $\beta > 0$.}$$
\end{itemize}

\noindent Travel time evaluation is a straightforward calculation of the expected travel time.
More interestingly, perceived travel time evaluation combines travel time with a \emph{dispersion} term that is weighted by the parameter $\beta > 0$.
This evaluation function was inspired by \citet{anderson1988representative}, who introduce a similar function and show that it has favorable properties in the logit context.
The dispersion term captures a distributional property of the probability vector $p \in \mathcal{P}$.
As such, perceived travel time favors probability vectors that are spread out over the routes, which creates robustness against uncertainty.

\paragraph*{Route Set Measures}
Finally, we obtain a \emph{measure} of route set quality by combining a routing (a model for traveler behavior) with an evaluation function (a scoring function for given behavior).
That is, for routing $\routingBold(\routes)$ and evaluation function $\mathcal{E}(\routes, p)$, we define the measure $\mathcal{M}(\routes)=\mathcal{E}(\routes, \routingBold(\routes))$.
Table~\ref{tbl:measures} summarizes the four measures that can be obtained by substituting one of the two routings into one of the two evaluation functions.
For shortest path routing with travel time evaluation, the measure is simply the duration $l_{\min}$ of the shortest route.
The same measure is obtained when perceived travel time is used, because the dispersion term evaluates to zero when the full probability is assigned to a single route.\footnote{Using the convention that $p_i \log(p_i) = 0$ for $p_i=0$.}
The measure for logit routing with travel time evaluation corresponds to calculating the expected duration with logit probabilities.
Logit routing with perceived travel time simplifies to a logsum that is scaled by a factor $-\frac{1}{\beta}$.

\newcommand\measuretable{
	\begin{adjustbox}{width=\linewidth,center}
		\begin{tabular}{lcccccc}
			\toprule
			& \multicolumn{3}{c}{Shortest Path Routing} & \multicolumn{3}{c}{Logit Routing}\\
			\cmidrule(lr){2-4} \cmidrule(lr){5-7}
			Evaluation & Measure & Monot. & Consist. & Measure & Monot. & Consist.\\
			\midrule
			Travel Time & $l_\text{min}$ & \cmark & \cmark & $\sum_{i\in \routes} p^{logit}_i l_i$ & \xmark & \xmark\\
			Perceived Travel Time & $l_\text{min}$ & \cmark & \xmark & $-\frac{1}{\beta}\log \left(\sum_{i\in \routes} e^{-\beta l_i}\right)$ & \cmark \cmark & \cmark\\
			\bottomrule
		\end{tabular}
	\end{adjustbox}
	\caption{Overview of route set measures and properties (double checkmark for strict monotonic).}
}
\begin{table}[t]
	\centering
	\measuretable{}
	\label{tbl:measures}
\end{table}

\begin{figure}[t]
	\centering
	\includegraphics[scale=0.5]{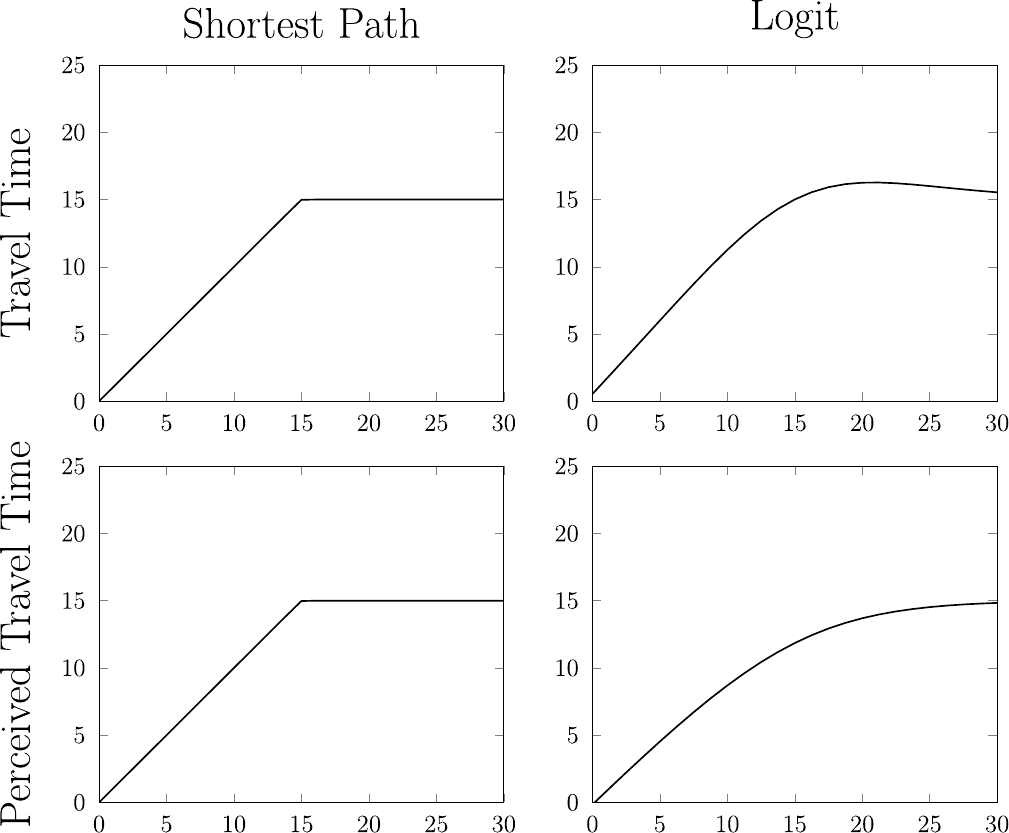}
	\caption{Example measures (y-axis) for $\routes = \{1,2\}$, $l_1=15$, $\beta=0.22$, and varying $l_2$ (x-axis).}
	\label{fig:measures_plots}
\end{figure}

While logit routing with travel time evaluation has been used in the literature (e.g., \citet{bertsimas2020joint,HartlebSchmidt2022-RailwayTimetablingIntegrated}), the example in Figure~\ref{fig:measures_plots} reveals some unexpected behavior:
When the duration of route 2 is increased from 20 minutes to 30 minutes, \emph{the measure improves while the route set has become worse}.
This shortcoming of logit routing with travel time evaluation is formalized in the next section.
The other route set measures are assigned a symbol, and will be used later as building blocks for the evaluation of timetables and line plans:
\begin{itemize}
	\item \emph{Route Set Measure for Shortest Path Routing with Travel Time Evaluation:} \begin{equation}\mathcal{M}^{RS}_{sp}(\routes) = l_{\min}.\label{eq:M_RS_sp}\end{equation}
	\item \emph{Route Set Measure for Logit Routing with Perceived Travel Time Evaluation:} \begin{equation}\mathcal{M}^{RS}_{logit}(\routes) = -\frac{1}{\beta}\log \left(\sum_{i\in \routes} e^{-\beta l_i}\right).\label{eq:M_RS_logit}\end{equation}
\end{itemize}

\paragraph*{Monotonicity and Consistency}
This section introduces two desirable properties of route set measures: monotonicity and consistency.
First, we define a partial order on route sets to capture the intuition that a route set is better when the routes are faster and more options are available.
In particular, consider two route sets $\routes$ and $\routes'$ and corresponding route duration vectors $l$ and $l'$.
Now sort $l$ and $l'$ and append routes with infinite duration to make the vectors the same size.
We say that $\routes \preceq \routes'$ if $l \le l'$ elementwise, and $\routes \prec \routes'$ if additionally $l \neq l'$.
Informally, route set $\routes$ is better when each route in $\routes'$ can be matched to a better route in $\routes$.

The first desirable property is monotonicity.
We say that a route set measure is \emph{monotonic} if better route sets (according to the partial order defined above) receive better scores.
The second desirable property is consistency.
A route set measure is \emph{consistent} when the travelers behave in such a way that they optimize the evaluation function.
This suggests that the route set is evaluated in a way that matches what is important to the travelers.
When a measure is inconsistent, the performance is measured according to a metric that the travelers do not agree with, as revealed through their routing choices.
These properties are formalized as follows:
\begin{itemize}
	\item \emph{Monotonicity:}
	A route set measure is monotonic if, for all route sets $\routes$ and $\routes'$,
	\begin{equation*}
		\routes \preceq \routes' \Longrightarrow \mathcal{M}(\routes) \le \mathcal{M}(\routes'),
	\end{equation*}
	and strict monotonic if additionally
	\begin{equation*}
		\routes \prec \routes' \Longrightarrow \mathcal{M}(\routes) < \mathcal{M}(\routes').
	\end{equation*}
	\item \emph{Consistency:}
	A route set measure that combines routing $\routingBold(\routes)$ and evaluation function $\mathcal{E}(\routes, p)$ is \emph{consistent} if, for any route set $\routes$,
	\begin{equation*}
		\routingBold(\routes) \textrm{ is an optimal solution to the minimization problem } \min_{p \in \mathcal{P}} \mathcal{E}(\routes, p).
	\end{equation*}
\end{itemize}

Table~\ref{tbl:measures} summarizes the monotonicity and consistency properties of the four route set measures, with proofs provided in Appendix~\ref{app:properties_proofs}.
The table shows that shortest path routing with travel time evaluation and logit routing with perceived travel time evaluation both exhibit monotonicity and consistency.
The counter-intuitive results for logit routing with travel time evaluation are explained by a lack of monotonicity, and the measure lacks consistency as well.
Finally, we remark that monotonicity alone is not sufficient to obtain a good measure: for example, $\mathcal{M}(\routes)=0$ is monotonic but practically useless.
Consistency is also useful to motivate the assumptions that may appear in other parts of the analysis.
For example, Table~\ref{tbl:measures} shows that shortest path routing with travel time evaluation and shortest path routing with perceived travel time evaluation have the same measure, but using travel time yields a more consistent theory.

\paragraph*{Random Utility Model Interpretation}
\label{sec:ptt_motivation}
The logit model is often motivated by formulating it as a random utility model \citep{Train2009-DiscreteChoiceMethods}.
This perspective is also useful in the current context, and we introduce this model here as the \textit{random travel time model}.
In the random travel time model, the \textit{perceived duration} $\chi_i$ of a route is given by
\begin{equation*}
	\chi_i = l_i - \frac{1}{\beta}\varepsilon_i, \quad \forall i \in \routes,
\end{equation*}
with independent error terms $\varepsilon_i$ that are distributed according to a Gumbel distribution.
Each traveler chooses the route that minimizes perceived duration \emph{for them}.
This choice is motivated by the route durations $l_i$, but also by a personal random component $\varepsilon_i$ for each route.
Again, for $\beta\rightarrow\infty$ the randomness disappears and travelers will pick the shortest path.

Analogous to the random utility model, for the random travel time model, it can be shown that
\begin{equation*}
	\mathbb{E}\left(\min_{i \in \routes} \chi_i\right) = -\frac{1}{\beta}\log \left(\sum_{i\in \routes} e^{-\beta l_i}\right) + C,
\end{equation*}
where $C$ is an unknown constant that reflects the base utility of the travelers \citep{Train2009-DiscreteChoiceMethods}.
Note that (up to a constant) this quantity is precisely $\mathcal{M}^{RS}_{logit}$.
This shows that the measure that combines perceived travel time and logit routing can be interpreted as the expected perceived travel time of a random traveler.

The term $\log \left(\sum_{i\in \routes} e^{-\beta l_i}\right)$ is referred to as the \emph{logsum} in the literature, and has been used extensively to assess the user benefits of public transport. However, as discussed in Section~\ref{sec:lit}, existing approaches compute the logsum change compared to a baseline scenario, which can subsequently be converted to units of time. Our approach is to directly convert the logsum to units of time, which eases its interpretation and omits the need for a baseline scenario. 
\subsection{Timetables}
\label{sec:tt}

We now shift attention to evaluating timetables for systems that provide periodic service.
Such a system is defined by a route set $\routes=\{1,\hdots,n\}$, a period $T > 0$, and a timetable $\timetableBold \in \left[0,T\right)^n$ of departure times (see Figure~\ref{fig:concepts}).
The routes are performed periodically, such that route $i \in \routes$ departs at times $\timetable_i$, $\timetable_i + T$, $\timetable_i + 2T$, etc.
Rather than a static route set, travelers now face an \emph{observed route set} with added waiting or \textit{adaptation} times that depend on the preferred departure —or, equivalently, preferred arrival— time of the travelers, which we also refer to as \textit{travel demand}. Since travelers typically plan their journey according to the timetable, waiting is only partially visible at stations. However, unless a traveler's preferred departure time exactly coincides with a timetabled departure time, waiting is inevitable, whether at home, in the office, or at the station. 
We will evaluate the experience \emph{for a particular traveler} by evaluating the observed route set with the measures developed in the previous section. 
To evaluate the full timetable, we propose to calculate the \emph{expected observed route set measure} for randomly distributed travel demand, reflecting the average experience of travelers facing this timetable.

This section first introduces the observed route set measure, and studies how it changes as a function of the preferred departure time.
We show how to calculate characteristic values that describe this function, and how to use them to calculate timetable measures under uniformly distributed travel demand.
These measures are further specialized to shortest path routing and logit routing, respectively, to obtain measures that are completely separable and convex in the characteristic values, paving the way for extensions to line planning.

The following notation is used to work with cyclic timetables.
The timetable $\timetableBold \in \left[0,T\right)^n$ defines a natural cyclic ordering of the routes by departure time.
For a route $i \in \routes$, let $\pi(i)$ be its predecessor and let $\sigma(i)$ be its successor.
Routes that depart at the same time are ordered according to their index.
For example, if subsequent routes $i$ and $\sigma(i)$ depart at the same time $\timetable_i = \timetable_{\sigma(i)}$, then it must be that $i < \sigma(i)$.
Time intervals $[t,t']$ are also considered to be cyclical.
That is, $[t,t']$ equals $[t,t']$ when $t \le t'$ and $[t,T) \cup [0, t']$ otherwise. For example, if $T=60$, the interval $[55,5]$ equals $[55,60) \cup [0, 5]$.

\paragraph*{Observed Route Set Measures}
In the periodic setting, a traveler wishing to depart at time $t$ has to take waiting time into account.
The waiting time for route $i \in \routes$ is given by $\timetable_i + zT - t$ for the smallest value of $z \in \mathbb{Z}$ such that the waiting time is non-negative.
For convenience, we define the function $[x]_T$ to be the modulo-like operation that maps $x$ onto $[0,T)$ by adding or subtracting a multiple of $T$.
Waiting time for route $i \in \routes$ is then simply $[\timetable_i-t]_T$.
We define the \emph{observed route set} $\routes_\timetable(t)$ to be the route set as experienced by a traveler arriving at time $t$.
That is, the observed route durations are given by $l_i + [\timetable_i - t]_T$ for every route $i \in \routes$.

We are interested in studying the properties of the \emph{observed route set measure} $\mathcal{M}^{RS}(\routes_\timetable(t))$ as a function of time, where $\mathcal{M}^{RS}$ represents either the route set measure for shortest path \eqref{eq:M_RS_sp} or for logit \eqref{eq:M_RS_logit}.
That is, we consider the following two functions:
\begin{align}
	\mathcal{M}^{RS}_{sp}(\routes_\timetable(t)) & = \min_{i\in \routes}\left\{l_i + [\timetable_i - t]_T\right\}, \label{eq:rs_sp}\\
	\mathcal{M}^{RS}_{logit}(\routes_\timetable(t)) & = -\frac{1}{\beta}\log \left(\sum_{i\in \routes} e^{-\beta \left(l_i + [\timetable_i - t]_T\right)}\right).\label{eq:rs_logit}
\end{align}
The reason to study these functions is to later derive simple expressions for the expected observed route set measure $\mathbb{E}_t\left(\mathcal{M}^{RS}(\routes_\timetable(t))\right)$, which will be used as a measure for timetable quality.

Figure~\ref{fig:timetable_measure_example} visualizes the observed route set measure for shortest path routing.
Figure~\ref{fig:timetable_measure_example1} shows an hourly timetable, where circled numbers indicate departures of a given duration that are located according to the schedule. 
For example, the number 10 indicates that a route of duration 10 departs at time 50, so every hour at xx:50. 
Figure~\ref{fig:timetable_measure_example2} plots the function $\mathcal{M}^{RS}_{sp}(\routes_\timetable(t))$.
At $t=\theta_1=5$ the shortest path is to use route 1, which departs immediately.
The corresponding duration is $l_1=20$ and the waiting time is zero, which results in an observed route set measure of $\mathcal{M}^{RS}_{sp}(\routes_\timetable(5)) = 20$.
The choice for route 1 is optimal everywhere in the interval $(50, 5]$.
At $t=55$, for example, it is optimal to wait for 10 minutes and then take route 1, which results in an observed route set measure of $\mathcal{M}^{RS}_{sp}(\routes_\timetable(55)) = 30$.
Route 2 ($l_2=30$) is sub-optimal: It is always better to wait an additional 10 minutes to take route 3 ($l_3 = 15$).
This is reflected in the plot by the lack of a jump at $\timetable_2 = 10$, as the route choice does not change.

\begin{figure}[t]
\def\lOne{20}
\def\lTwo{30}
\def\lThree{15}
\def\lFour{10}
\def\thetaOne{5}
\def\thetaTwo{10}
\def\thetaThree{20}
\def\thetaFour{50}
\centering
\begin{subfigure}[b]{0.35\textwidth}
\centering
\begin{tikzpicture}[>=latex]
	\coordinate (O) at (0,0);
	\def\radius{1.7cm}
	\def\smallRadius{0.5cm}
	  
	\coordinate (left) at (180:\smallRadius);
	\coordinate (right) at (0:\smallRadius); 
	\coordinate (top) at (90:\smallRadius);
	\coordinate (bottom) at (-90:\smallRadius);
	\draw[->,thick,out=90,in=180] (left) to (top);
	\draw[->,thick,out=0,in=90] (top) to (right);
	\draw[->,thick,out=-90,in=0] (right) to (bottom);
	\draw[->,thick,out=180,in=-90] (bottom) to (left);
	
	\draw (O) circle[radius=\radius];
	
	\coordinate (0t) at (90:\radius); 
	\fill (0t)  (90: \radius-0.32cm) node{$0$};
	\coordinate (180t) at (-90:\radius); 
	\fill (180t)  (-90: \radius-0.35cm) node{30};
	\coordinate (15t) at (0:\radius); 
	\fill (15t)  (0: \radius-0.39cm) node{15};
	\coordinate (45t) at (180:\radius); 
	\fill (45t)  (180: \radius-0.39cm) node{45};
	
	\coordinate (e1) at (90-\thetaOne*6:\radius);
	\coordinate (e2) at (90-\thetaTwo*6:\radius);
	\coordinate (e3) at (90-\thetaThree*6:\radius);
	\coordinate (e4) at (90-\thetaFour*6:\radius);
	
	\node[vertex, fill=white] at (e1) {\lOne};
	\node[vertex, fill=white] at (e2) {\lTwo};
	\node[vertex, fill=white] at (e3) {\lThree};
	\node[vertex, fill=white] at (e4) {\lFour};
	
	\node at (0,-2.3) {\textcolor{white}{text}};
\end{tikzpicture}
\caption{Timetable.}
\label{fig:timetable_measure_example1}
\end{subfigure}
\hspace{0.01\textwidth}
\begin{subfigure}[b]{0.50\textwidth}
\centering
\begin{tikzpicture}[x=5cm,y=0.1cm]
	\draw [thick] (0,40) -- (0,0) -- (1,0);
	\foreach \x in {0, 10, 20, 30, 40}
	\draw[thick] (-0.5/40, \x) -- (0, \x);
	\foreach \x in {5, 10, 15, 20, 25, 30, 35, 40}
	\draw[thin,dotted] (0,\x) -- (1,\x);  
	\foreach \x/\xtext in {10, 20, 30, 40}
	\node [left] at (0,\x) {$\xtext$};
	\foreach \x in {0, 15, 30, 45, 60}
	\draw[thick] (\x/60, 0) -- (\x/60, -0.5);
	\foreach \x in {0, 5, 10, 15, 20, 25, 30, 35, 40, 45, 50, 55}
	\draw[thin,dotted] (\x/60, 40) -- (\x/60, 0);
	
	\node [left] at (0,-2) {0};
	
	\node [left] at (-0.15,18) {$\mathcal{M}^{RS}_{sp}$};
	\node [below] at (0.25,0) {15};
	\node [below] at (0.5,0) {30};
	\node [below] at (0.75,0) {45};
	\node [below] at (1,0) {60};
	\node [below] at (0.5,-4) {$t$};
	
	\coordinate (tau1) at (\thetaOne/60,\lOne);
	\coordinate (tau2) at (\thetaTwo/60,\lTwo);
	\coordinate (tau3) at (\thetaThree/60,\lThree);
	\coordinate (tau4) at (\thetaFour/60,\lFour);
	
	\draw[thick] (0,\lOne+\thetaOne-0) -- (tau1);
	\draw[fill=black] (tau1) circle (2pt);
	\draw[thick] (\thetaOne/60,\lThree+\thetaThree-\thetaOne) -- (\thetaTwo/60,\lThree+\thetaThree-\thetaTwo); 
	\draw[fill=white] (\thetaOne/60,\lThree+\thetaThree-\thetaOne) circle (2pt);
	\draw[fill=black] (\thetaTwo/60,\lThree+\thetaThree-\thetaTwo) circle (2pt);
	\draw[thick] (\thetaTwo/60,\lThree+\thetaThree-\thetaTwo) -- (tau3); 
	\draw[fill=black] (tau3) circle (2pt);
	\draw[thick] (\thetaThree/60,\lFour+\thetaFour-\thetaThree) -- (tau4);	
	\draw[fill=white] (\thetaThree/60,\lFour+\thetaFour-\thetaThree) circle (2pt);
	\draw[fill=black] (tau4) circle (2pt);
	\draw[thick] (\thetaFour/60,\lOne+\thetaOne+60-\thetaFour) -- (60/60,\lOne+\thetaOne);
	\draw[fill=white] (\thetaFour/60,\lOne+\thetaOne+60-\thetaFour) circle (2pt);

	
\end{tikzpicture}
\caption{Observed route set measure $\mathcal{M}^{RS}_{sp}(\routes_\timetable(t))$.}
\label{fig:timetable_measure_example2}
\end{subfigure}

\caption{Example timetable and observed shortest path route set measure for route set $\routes=\{1,2,3,4\}$, route lengths $l=(20,30,15,10)$, period $T=60$, and schedule $\timetableBold=(5,10,20,50)$.}
\label{fig:timetable_measure_example}
\end{figure}

The example also displays two useful properties that are true in general, as proven in Appendix~\ref{app:tt_proofs}.\ref{prop:constantrouting}-\ref{prop:trans_invar}:
\begin{itemize}
	\item \emph{Constant Routing between Departures:}
	For shortest path routing and logit routing the route choice does not change between departures.
	That is, for route $i \in \routes$ and successor $\sigma(i) \in \routes$ it holds that $\routingBold(\routes_\timetable(t))$ is constant for $t \in (\timetable_i, \timetable_{\sigma(i)}]$.
	\item \emph{Translation Invariance:}
	Both $\mathcal{M}^{RS}_{sp}$ and $\mathcal{M}^{RS}_{logit}$ are translation invariant between departures.
	That is, for route $i \in \routes$ and successor $\sigma(i) \in \routes$ it holds that $$\mathcal{M}^{RS}(\routes_\timetable(t)) = [\theta_{\sigma(i)} - t]_T + \mathcal{M}^{RS}(\routes_\timetable(\timetable_{\sigma(i)})), \quad \forall t \in (\timetable_i, \timetable_{\sigma(i)}].$$
\end{itemize}
The first property informally states that travelers do not make routing decisions at the time they prefer to depart, but postpone their decisions until the first route departs, essentially treating the waiting time until the first departure as sunk cost.
As a direct consequence, the second property states that the observed route set measure can be split into two components: the waiting time until the first departure and the observed route set measure at the first departure. 
During the waiting part, the slope of the observed route set measure is $-1$, as can also be seen in Figure~\ref{fig:timetable_measure_example2}.

\paragraph*{Characteristic Values}
Given that the observed route set measure has a fixed slope between departures, we can represent $\mathcal{M}^{RS}(\routes_\timetable(t))$ with a finite number of values.
We define the following characteristic values:
\begin{itemize}
	\item $\delta_i$: Waiting time between predecessor route $\pi(i)$ and route $i \in \routes$.
	\item $\tau_i$: Observed route set measure at time $\timetable_i$ when $i \in \routes$ is the next departure.
	\item $\Delta_i$: Change in observed route set measure due to missing route $i \in \routes$.
\end{itemize}
The characteristic values are visualized by Figure~\ref{fig:representations}.
There are two ways to fully describe $\mathcal{M}^{RS}(\routes_\timetable(t))$ with the characteristic values, up to shifts of the complete schedule.
The $\delta$-representation (Figure~\ref{fig:delta_repres}) defines the function based on the gaps $\delta_i$ between departure times, and the observed route set measures $\tau_i$ at these departures times.
Using instead the vertical differences, the $\Delta$-representation (Figure~\ref{fig:Delta_repres}) describes the function with the observed route set measures $\tau_i$ at each departure time, and the change $\Delta_i$ in value that would result from missing this route. Note that since route 2 is sub-optimal, $\Delta_2=0$ in  Figure~\ref{fig:Delta_repres}.

\begin{figure}[t]
\def\lOne{20}
\def\lTwo{30}
\def\lThree{15}
\def\lFour{10}
\def\thetaOne{5}
\def\thetaTwo{10}
\def\thetaThree{20}
\def\thetaFour{50}
\centering
\begin{subfigure}[b]{0.47\textwidth}
\centering
\begin{tikzpicture}[x=5cm,y=0.1cm]
	\draw [thick] (0,40) -- (0,0) -- (1,0);
	\foreach \x in {0, 10, 20, 30, 40}
	\draw[thick] (-0.5/40, \x) -- (0, \x);
	\foreach \x/\xtext in {10, 20, 30, 40}
	\node [left] at (0,\x) {$\xtext$};
	\foreach \x in {0, \thetaOne, \thetaTwo, \thetaThree, \thetaFour, 60}
	\draw[thick] (\x/60, 0) -- (\x/60, -0.5);
	
	\node [left] at (0,-2) {0};
	\node [left] at (-0.12,18) {$\mathcal{M}^{RS}_{sp}$};
	\node [below] at (\thetaOne/60, 0){$\theta_1$};
	\node [below] at (\thetaTwo/60, 0){$\theta_2$};
	\node [below] at (\thetaThree/60, 0){$\theta_3$};
	\node [below] at (\thetaFour/60, 0){$\theta_4$};
	\node [below] at (0.5,-4) {$t$};
	
	\coordinate (tau1) at (\thetaOne/60,\lOne);
	\coordinate (tau2) at (\thetaTwo/60,\lTwo);
	\coordinate (tau3) at (\thetaThree/60,\lThree);
	\coordinate (tau4) at (\thetaFour/60,\lFour);
	
	\draw[thick] (0,\lOne+\thetaOne-0) -- (tau1);
	\draw[fill=black] (tau1) circle (2pt);
	\draw[thick] (\thetaOne/60,\lThree+\thetaThree-\thetaOne) -- (\thetaTwo/60,\lThree+\thetaThree-\thetaTwo); 
	\draw[fill=white] (\thetaOne/60,\lThree+\thetaThree-\thetaOne) circle (2pt);
	\draw[fill=black] (\thetaTwo/60,\lThree+\thetaThree-\thetaTwo) circle (2pt);
	\draw[thick] (\thetaTwo/60,\lThree+\thetaThree-\thetaTwo) -- (tau3); 
	\draw[fill=black] (tau3) circle (2pt);
	\draw[thick] (\thetaThree/60,\lFour+\thetaFour-\thetaThree) -- (tau4);	
	\draw[fill=white] (\thetaThree/60,\lFour+\thetaFour-\thetaThree) circle (2pt);
	\draw[fill=black] (tau4) circle (2pt);
	\draw[thick] (\thetaFour/60,\lOne+\thetaOne+60-\thetaFour) -- (60/60,\lOne+\thetaOne);
	\draw[fill=white] (\thetaFour/60,\lOne+\thetaOne+60-\thetaFour) circle (2pt);
	
	\tikzstyle{redarrow}=[red, latex-latex, thick]
	\tikzstyle{redarrowlr}=[red, -latex, thick]
	\draw[redarrowlr] (0/60, 20) -- (\thetaOne/60, 20);
	\draw[redarrowlr] (60/60, 20) -- node[below] {$\delta_1$} (\thetaFour/60, 20);
	\draw[redarrow] (\thetaOne/60, 25) -- node[below] {$\delta_2$} (\thetaTwo/60, 25);
	\draw[redarrow] (\thetaTwo/60, 15) -- node[below] {$\delta_3$} (\thetaThree/60, 15);
	\draw[redarrow] (\thetaThree/60, 10) -- node[below] {$\delta_4$} (\thetaFour/60, 10);
	\draw[redarrow] (\thetaThree/60, 10) -- node[right] {$\delta_4$} (\thetaThree/60, 40);
	\draw[redarrow] (\thetaFour/60, 0) -- node[right] {$\tau_4$} (\thetaFour/60, 10);
	
	\path [fill=red,fill opacity=0.1] (\thetaThree/60, 0) -- (\thetaThree/60, 40) -- (\thetaFour/60, 10) -- (\thetaFour/60, 0);
	
	
\end{tikzpicture}
\caption{$\delta$-representation.}
\label{fig:delta_repres}
\end{subfigure}
\hspace{0.01\textwidth}
\begin{subfigure}[b]{0.47\textwidth}
\centering
\begin{tikzpicture}[x=5cm,y=0.1cm]
	\draw [thick] (0,40) -- (0,0) -- (1,0);
	\foreach \x in {0, 10, 20, 30, 40}
	\draw[thick] (-0.5/40, \x) -- (0, \x);
	\foreach \x/\xtext in {10, 20, 30, 40}
	\node [left] at (0,\x) {$\xtext$};
	\foreach \x in {0, \thetaOne, \thetaTwo, \thetaThree, \thetaFour, 60}
	\draw[thick] (\x/60, 0) -- (\x/60, -0.5);
	\node [left] at (-0.12,18) {$\mathcal{M}^{RS}_{sp}$};
	\node [left] at (0,-2) {0};
	
	\node [below] at (\thetaOne/60, 0){$\theta_1$};
	\node [below] at (\thetaTwo/60, 0){$\theta_2$};
	\node [below] at (\thetaThree/60, 0){$\theta_3$};
	\node [below] at (\thetaFour/60, 0){$\theta_4$};
	\node [below] at (0.5,-4) {$t$};
	
	\coordinate (tau1) at (\thetaOne/60,\lOne);
	\coordinate (tau2) at (\thetaTwo/60,\lTwo);
	\coordinate (tau3) at (\thetaThree/60,\lThree);
	\coordinate (tau4) at (\thetaFour/60,\lFour);
	
	\draw[thick] (0,\lOne+\thetaOne-0) -- (tau1);
	\draw[fill=black] (tau1) circle (2pt);
	\draw[thick] (\thetaOne/60,\lThree+\thetaThree-\thetaOne) -- (\thetaTwo/60,\lThree+\thetaThree-\thetaTwo); 
	\draw[fill=white] (\thetaOne/60,\lThree+\thetaThree-\thetaOne) circle (2pt);
	\draw[fill=black] (\thetaTwo/60,\lThree+\thetaThree-\thetaTwo) circle (2pt);
	\draw[thick] (\thetaTwo/60,\lThree+\thetaThree-\thetaTwo) -- (tau3); 
	\draw[fill=black] (tau3) circle (2pt);
	\draw[thick] (\thetaThree/60,\lFour+\thetaFour-\thetaThree) -- (tau4);	
	\draw[fill=white] (\thetaThree/60,\lFour+\thetaFour-\thetaThree) circle (2pt);
	\draw[fill=black] (tau4) circle (2pt);
	\draw[thick] (\thetaFour/60,\lOne+\thetaOne+60-\thetaFour) -- (60/60,\lOne+\thetaOne);
	\draw[fill=white] (\thetaFour/60,\lOne+\thetaOne+60-\thetaFour) circle (2pt);
	
	\tikzstyle{redarrow}=[red, latex-latex, thick]
	\tikzstyle{redarrowlr}=[red, -latex, thick]
	\draw[redarrow]  (\thetaOne/60, 20) node[right]{$\Delta_1$} -- (\thetaOne/60, 30);
	\draw[fill=red]  (\thetaTwo/60, 25) circle(2pt) node[right, red]{$\Delta_2$};
	\draw[redarrow]  (\thetaThree/60, 15) -- node[right]{$\Delta_3$} (\thetaThree/60, 40);
	\draw[redarrow]  (\thetaFour/60, 10) -- node[right]{$\Delta_4$} (\thetaFour/60, 35);
	\draw[redarrow] (\thetaThree/60, 0) -- node[right] {$\tau_3$} (\thetaThree/60, 15);
	\draw[redarrow] (\thetaFour/60, 0) -- node[right] {$\tau_4$} (\thetaFour/60, 10);

	
\end{tikzpicture}
\caption{$\Delta$-representation.}
\label{fig:Delta_repres}
\end{subfigure}

\caption{Two representations for the observed route set measure $\mathcal{M}^{RS}(\routes_\timetable(t))$ in Figure~\ref{fig:timetable_measure_example2}.}
\label{fig:representations}
\end{figure}

Algorithm~\ref{alg:characteristic} provides pseudocode to calculate the characteristic values.
The calculation of $\delta$ (Line~\ref{line:delta}) and $\tau$ (Lines~\ref{line:tau_sp} and \ref{line:tau_logit}) follows immediately from the definitions, except for a waiting time correction in the case of simultaneous departures (Line~\ref{line:correction}).
This technical detail is justified in Appendix~\ref{app:tt_proofs}.\ref{prop:alg_char}.
The value of the measure after missing route $i \in \routes$ is conveniently expressed as $\delta_{\sigma(i)} + \tau_{\sigma(i)}$ by virtue of translation invariance (visualized in Figure~\ref{fig:delta_repres} for $\delta_4 + \tau_4$).
It follows that $\Delta_i$ can be calculated as $\delta_{\sigma(i)} + \tau_{\sigma(i)} - \tau_i$ (Line~\ref{line:Delta}).

\begin{figure}[t]
\centering
\begin{algorithm}[H]
\caption{Calculating the Characteristic Values}\label{alg:characteristic}
	\begin{algorithmic}[1]
		\Function{$w$}{$i,j$} \Comment{Waiting time correction for simultaneous departures} \label{line:correction}
			\If{$\theta_i = \theta_j$ and $j < i$}
				\State \textbf{return} $T$ \Comment{Add full period $T$ if $j$ is processed before $i$}
			\Else
				\State \textbf{return} $[\theta_j - \theta_i]_T$ \Comment{Regular wait time}
			\EndIf
		\EndFunction\\
		\Procedure{CharacteristicValues}{$\routes$, $T$, $\timetable$}
			\If{$\lvert \routes \rvert = 1$}
				\State $\delta_1 \gets T$
				\State $\tau_1 \gets l_1$
			\Else
				\For{$i \in \routes$}
					\State $\delta_i \gets w(\pi(i), i)$ \Comment{From definition} \label{line:delta}
					\If{shortest path routing}
						\State $\tau_i \gets \min_{j \in \routes}\{l_j + w(i,j)\}$ \Comment{From route set measure~\eqref{eq:rs_sp}} \label{line:tau_sp}
					\ElsIf{logit routing}
						\State $\tau_i \gets -\frac{1}{\beta}\log\left(\sum_{j\in \routes} e^{-\beta (l_j + w(i,j))} \right)$ \Comment{From route set measure~\eqref{eq:rs_logit}} \label{line:tau_logit}
					\EndIf
				\EndFor
			\EndIf
			\For{$i \in \routes$}
				\State $\Delta_i \gets \delta_{\sigma(i)} + \tau_{\sigma(i)} - \tau_i$ \Comment{From translation invariance} \label{line:Delta}
			\EndFor
			\State \textbf{return} $\delta$, $\tau$, $\Delta$
		\EndProcedure
	\end{algorithmic}
\end{algorithm}
\end{figure}

\paragraph*{Timetable Measures}
To measure the quality of timetables, we propose to use the \emph{expected observed route set measure}, $\mathbb{E}_t\left(\mathcal{M}^{RS}(\routes_\timetable(t))\right)$.
This measure can be interpreted as the expected (perceived) travel time experienced by a random traveler, or as a weighted average over all travelers.
As a weighted average of route set measures, the timetable measure inherits (strict) monotonicity, and improves when the underlying routes are reduced in duration or when new routes are added.
The measure also remains consistent, in the sense that the timetable is evaluated in a way that now matches what is important to the \emph{average} traveler.
\begin{itemize}
	\item \emph{Timetable Measure for Shortest Path Routing or Logit Routing (general demand):}
	\begin{equation}
		\mathcal{M}^{TT}(\routes, T, \timetableBold) = \mathbb{E}_t\left(\mathcal{M}^{RS}(\routes_\timetable(t))\right).
	\end{equation}
\end{itemize}
Note that the function $\mathcal{M}^{RS}(\routes_\timetable(t))$ has been completely specified previously.
As a result, the timetable measure can easily be estimated through sampling or numerical integration.

Furthermore, the measure can be expressed in closed form if we assume that travel demand is distributed uniformly over the period.
This is a common assumption in periodic timetabling and fairly accurate when the period $T$ is reasonably small (e.g., see \cite{gentile2016modelling,kaspi2013service,polinder2021timetabling}). Moreover, this assumption naturally reflects the idea that timetables with constant headways are preferred over those with uneven headways, as they reduce expected waiting time. Using the characteristic values calculated with Algorithm~\ref{alg:characteristic}, it is proven in Appendix~\ref{app:tt_proofs}.\ref{prop:tt_unif_general} that the timetable measure can be calculated in two distinct ways that correspond to the $\delta$- and $\Delta$-representation, respectively:

\begin{itemize}
	\item \emph{Timetable Measure for Shortest Path Routing or Logit Routing (uniform demand):}
	\begin{equation}
		\label{eq:tt_uniform}
		\mathcal{M}^{TT}(\routes, T, \timetableBold) = \frac{1}{T} \sum_{i\in \routes} \left( \frac{1}{2} \delta_i^2 + \tau_i \delta_i \right) = \frac{1}{T} \sum_{i\in \routes} \left( \frac{1}{2} \Delta_i^2 + \tau_i \Delta_i \right).
	\end{equation}
\end{itemize}

Finally, we provide expressions for the timetable measures that are completely separable and convex.
These properties are extremely important when timetable measures are used as a building block, and they enable the efficient calculation of line plan measures in the next section.
For shortest path routing, note the following: if for route $i \in \routes$ we have that $\mathcal{M}^{RS}_{sp}(\routes_\timetable(\timetable_i)) < l_i$, then it means that route $i$ is not even the best option when it is about to depart, and the route can safely be removed without affecting the measure.
After removing the suboptimal routes, we have $\mathcal{M}^{RS}_{sp}(\routes_\timetable(\timetable_i)) = l_i$ for all routes, such that $\tau_i$ in \eqref{eq:tt_uniform} may be substituted for $l_i$.
For logit routing, $\tau_i$ may be replaced by an expression in $\Delta_i$ that is proven in Appendix~\ref{app:tt_proofs}.\ref{prop:tau_Delta_def}.
We then obtain the following two expressions:
\begin{itemize}
	\item \emph{Timetable Measure for Shortest Path Routing after preprocessing (uniform demand):}
	\begin{equation}
		\label{eq:tt_sp_convex}
		\mathcal{M}^{TT}_{sp}(\routes, T, \timetableBold) = \frac{1}{T} \sum_{i\in \routes} \left( \frac{1}{2} \delta_i^2 + l_i \delta_i \right).
	\end{equation}
	\item \emph{Timetable Measure for Logit Routing (uniform demand):}
	\begin{equation}
		\label{eq:tt_logit_convex}
		\mathcal{M}^{TT}_{logit}(\routes, T, \timetableBold) = \frac{1}{T} \sum_{i\in \routes} \left( \frac{1}{2} \Delta_i^2 + l_i \Delta_i + \frac{\Delta_i}{\beta}\log\left(\frac{1 - e^{-\beta \Delta_i}}{1 - e^{-\beta T}}\right) \right).
	\end{equation}
\end{itemize}
The convexity of \eqref{eq:tt_sp_convex} in $\delta$ is trivial, and the convexity of \eqref{eq:tt_logit_convex} in $\Delta$ is proven in Appendix~\ref{app:tt_proofs}.\ref{prop:measure_logit_convex}.
The fact that \eqref{eq:tt_logit_convex} is convex demonstrates the importance of viewing the observed route set measure through the $\Delta$-representation.
In fact, it can be shown that reformulating $\tau_i$ using the $\delta$-representation leads to an expression that is \emph{not} generally convex.

\paragraph*{Route Choice Models for Timetables}
The proposed measures naturally correspond to route choice models for which travelers distribute according to the expected route choice of a random traveler.
Formally, the routing $\routingBold^{TT}(\routes, T, \timetableBold)$ for a timetable can be defined as follows:

\begin{itemize}
    \item \emph{Route Choice Model for Shortest Path Routing or Logit Routing (general demand):}
	\begin{equation}
		\routingBold^{TT}(\routes, T, \timetableBold) = \mathbb{E}_t\left(\routingBold(\routes_\timetable(t))\right).
	\end{equation}
\end{itemize}

For shortest path routing, it holds that after removing the suboptimal routes, random travelers always choose the first departing route after their preferred departure time.
This gives a simple expression in $\delta$ for the route choice under uniform travel demand:
\begin{itemize}
	\item \emph{Route Choice Model for Shortest Path Routing after preprocessing (uniform demand):}
	\begin{equation}
		\label{eq:tt_routing_sp_unif}
		\routingBold^{TT}_{sp}(\routes, T, \timetableBold) = \left(\delta_1/T,...,\delta_n/T\right).
	\end{equation}
\end{itemize}
It is not obvious how to derive a similar expressions for logit routing when travelers face an arbitrary timetable.
But surprisingly, the next section will show that when random travelers face an \emph{optimal} timetable, they use route $i \in \routes$ with probability $\Delta_i/T$.
This again stresses the importance of the $\Delta$-representation.

\subsection{Line Plans}
\label{sec:lp}

Finally, we consider measures and route choice models for line plans.
A line plan is defined by a set of routes $\routes=\{1,\hdots,n\}$ that are performed every period $T>0$.
The exact timetable $\timetableBold \in \left[0,T\right)^n$ has not yet been decided, and will be optimized at a later stage (see Figure~\ref{fig:concepts}).

Because planners should be able to compare a large set of line plans, there is a need for line plan measures and route choice models that can be evaluated efficiently. Existing approaches typically rely on simplifying assumptions to achieve tractability \citep{Hartleb2022,schmidt2024planning}, including:
\begin{itemize}
    \item Ignoring waiting time entirely, which fails to capture the benefits of increased frequency (e.g., adding routes with the same travel time yields no improvement).
    \item Assuming equidistant departures, which implies an expected waiting time of $T/(2n)$. While easy to evaluate, this approach is sensitive to the ordering of routes and introduces non-monotonicity. 
\end{itemize} 

Because the timetable will eventually be optimized once the line plan is fixed, we propose to evaluate the quality of a line plan by calculating the minimum value of the timetable measure over all possible timetables. The corresponding routing is defined as the distribution of travelers under this optimal timetable.
When multiple OD pairs are involved, timetabling is known to be hard in theory and in practice \citep{liebchen2008first,lindner2022analysis}. This difficulty persists even in highly simplified settings: Appendix~\ref{app:lp_proofs}.\ref{prop:nphard} proves that when multiple OD pairs are involved, finding an optimal timetable is strongly NP-hard, even for a single origin and just two destinations.

To ensure tractability, the approach taken here is to not solve the full timetabling problem at once, but rather to optimize the timetable for each OD pair independently.
This section analyzes a single OD pair, but measures for multiple OD pairs can be combined, e.g., by taking a weighted average. The resulting measure represents the \textit{potential} of a line plan—an underestimation of the eventual (perceived) travel time once the timetable is finalized. As shown later in Section~\ref{subsec:lp-app}, this approximation performs remarkably well in practice on real-world networks.

To derive new measures and route choice models for line plans, we first discuss how to construct optimal timetables.
This requires solving different optimization problems, and it is shown how to do so efficiently.
The optimal solutions are then used to define the line plan measures and route choice models.

\paragraph*{Constructing Optimal Timetables}

The easiest way to construct an optimal timetable is to directly minimize the timetable measure for shortest path routing \eqref{eq:tt_sp_convex} or logit routing \eqref{eq:tt_logit_convex}.
Figure~\ref{fig:optimization_models} introduces two optimization models for this purpose.
The decision variables $x \ge 0$ represent the values of $\delta$ in \eqref{eq:tt_sp_convex}, and the decision variables $y \ge 0$ represent $\Delta$ in \eqref{eq:tt_logit_convex}.
The objectives \eqref{eq:sp:opt:obj} and \eqref{eq:logit:opt:obj} are chosen to match the corresponding timetable measure.
Finally, normalization constraints~\eqref{eq:sp:opt:total} and \eqref{eq:logit:opt:total} are added to exclude solutions that are not meaningful.
These constraints are justified by the fact that $\sum_{i \in \routes} \delta_i = T$ and $\sum_{i \in \routes} \Delta_i = T$ for any feasible timetable (Appendix~\ref{app:tt_proofs}.\ref{col:sum_T}).

\begin{figure}[t]
	\centering
	\begin{subfigure}{0.38\textwidth}
		\centering
		\begin{mini!}
			%
			{}
			%
			{\frac{1}{T} \sum_{i\in \routes} \left( \frac{1}{2} x_i^2 + l_i x_i \right), \label{eq:sp:opt:obj}}
			%
			{\label{formulation:sp}}
			%
			{}
			%
			\addConstraint
			{\sum_{i \in \routes} x_i}
			{= T, \label{eq:sp:opt:total}}
			{~}
			\addConstraint
			{x_i}
			{\ge 0 \quad \label{eq:sp:opt:xvar}}
			{\forall i \in \routes.}
		\end{mini!}%
		\caption{Shortest Path Routing.}
	\end{subfigure}%
	\hspace{0.01\textwidth}
	\begin{subfigure}{0.60\textwidth}
		\centering
		\begin{mini!}
			%
			{}
			%
			{\frac{1}{T} \sum_{i\in \routes} \left( \frac{1}{2} y_i^2 + l_i y_i + \frac{y_i}{\beta}\log\left(\frac{1 - e^{-\beta y_i}}{1 - e^{-\beta T}}\right) \right), \label{eq:logit:opt:obj}}
			%
			{\label{formulation:logit}}
			%
			{}
			%
			\addConstraint
			{\sum_{i \in \routes} y_i}
			{= T, \label{eq:logit:opt:total}}
			{~}
			\addConstraint
			{y_i}
			{\ge 0 \quad \label{eq:logit:opt:yvar}}
			{\forall i \in \routes.}
		\end{mini!}%
		\caption{Logit Routing.}
	\end{subfigure}
	\caption{Models for timetable optimization.}
	\label{fig:optimization_models}
\end{figure}

Given an optimal solution to Problem~\eqref{formulation:sp} or Problem~\eqref{formulation:logit}, Algorithm~\ref{alg:construct} provides methods to extend this information about $\delta$ or $\Delta$ into a complete timetable $\timetableBold$ that is feasible, and for which the measure matches the optimal objective value.
This is achieved by calculating desired values for the characteristic values, and then constructing the timetable accordingly.
It is not trivial that this is always possible, but Appendix~\ref{app:lp_proofs}.\ref{prop:match_sp}-\ref{prop:match_logit} proves that this is indeed the case.
Finally, Appendix~\ref{app:lp_proofs}.\ref{prop:optimal_sp}-\ref{prop:optimal_logit} proves that the output of Algorithm~\ref{alg:construct} is indeed an optimal timetable.
This is a consequence of the fact that the measure of the constructed timetable matches the optimal objective value.

\begin{figure}[t]
\centering
\begin{algorithm}[H]
\caption{Construct an Optimal Timetable}\label{alg:construct}
	\begin{algorithmic}[1]
		\Procedure{TimeTableShortestPath}{$x$} \Comment{Input: optimal solution to Problem~\eqref{formulation:sp}}
			\State $\timetable_1 \gets 0$ \Comment{Start timetable at arbitrary time 0}
			\For{$i \in \{2, \hdots, n\}$}
				\State $\hat{\delta}_i \gets x_i$ \Comment{Desired value for $\delta_i$ (from definition $x$)}
				\State $\timetable_i \gets \timetable_{\pi(i)} + \hat{\delta}_i$ \Comment{Progress timetable with desired gap between routes}
			\EndFor
			\State \textbf{return} $\timetable$
		\EndProcedure\\
		\Procedure{TimeTableLogit}{$y$} \Comment{Input: optimal solution to Problem~\eqref{formulation:logit}}
			\State $\timetable_1 \gets 0$ \Comment{Start timetable at arbitrary time 0}
			\For{$i \in \routes$}
				\State $\hat{\Delta}_i \gets y_i$ \Comment{Desired value for $\Delta_i$ (from definition $y$)}
				\State $\hat{\tau}_i \gets l_i + \frac{1}{\beta}\log\left(\frac{1 - e^{-\beta \hat{\Delta}_i}}{1 - e^{-\beta T}}\right)$ \Comment{Desired value for $\tau_i$ (from Appendix~\ref{app:tt_proofs}.\ref{prop:tau_Delta_def})}
			\EndFor
			\For{$i \in \{2, \hdots, n\}$}
				\State $\hat{\delta}_i \gets \hat{\tau}_{\pi(i)} + \hat{\Delta}_{\pi(i)} - \hat{\tau}_i$ \Comment{Desired value for $\delta_i$ (from translation invariance)}
				\State $\timetable_i \gets \timetable_{\pi(i)} + \hat{\delta}_i$ \Comment{Progress timetable with desired gap between routes}
			\EndFor
			\State \textbf{return} $\timetable$
		\EndProcedure
	\end{algorithmic}
\end{algorithm}
\end{figure}

\begin{figure}[t]
\centering
\begin{subfigure}[t]{\textwidth}
	\centering
	\begin{minipage}{0.35\textwidth}
		\centering
		\includegraphics{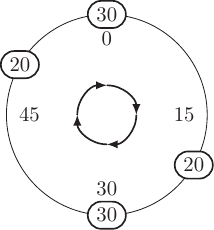}
	\end{minipage}
	\hspace{0.01\textwidth}
	\begin{minipage}{0.5\textwidth}
		\centering
		\def\lOne{30}
\def\lTwo{20}
\def\lThree{30}
\def\lFour{20}
\def\thetaOne{0}
\def\thetaTwo{20}
\def\thetaThree{30}
\def\thetaFour{50}
\begin{tikzpicture}[x=5cm,y=0.1cm]
	\draw [thick] (0,40) -- (0,10) -- (1,10);
	\foreach \x in {10, 20, 30}
	\draw[thick] (-0.5/40, \x) -- (0, \x);
	\foreach \x in {15, 20, 25, 30, 35, 40}
	\draw[thin,dotted] (0,\x) -- (1,\x);  
	\foreach \x/\xtext in {10, 20, 30, 40}
	\node [left] at (0,\x) {$\xtext$};
	\foreach \x in {0, 15, 30, 45, 60}
	\draw[thick] (\x/60, 10) -- (\x/60, 10-0.5);
	\foreach \x in {0, 5, 10, 15, 20, 25, 30, 35, 40, 45, 50, 55}
	\draw[thin,dotted] (\x/60, 40) -- (\x/60, 10);
	
	
	\node [left] at (-0.15,25) {$\mathcal{M}^{RS}_{sp}$};
	\node [below] at (0,10) {0};
	\node [below] at (0.25,10) {15};
	\node [below] at (0.5,10) {30};
	\node [below] at (0.75,10) {45};
	\node [below] at (1,10) {60};
	\node [below] at (0.5,10-4) {$t$};
	
	\coordinate (tau1) at (\thetaOne/60,\lOne);
	\coordinate (tau2) at (\thetaTwo/60,\lTwo);
	\coordinate (tau3) at (\thetaThree/60,\lThree);
	\coordinate (tau4) at (\thetaFour/60,\lFour);
	
	\draw[thick] (0,\lOne+\thetaOne-0) -- (tau1);
	\draw[fill=black] (tau1) circle (2pt);
	\draw[thick] (\thetaOne/60,\lTwo+\thetaTwo-\thetaOne) -- (tau2); 
	\draw[fill=white] (\thetaOne/60,\lTwo+\thetaTwo-\thetaOne) circle (2pt);
	\draw[fill=black] (tau2) circle (2pt);
	\draw[thick] (\thetaTwo/60,\lThree+\thetaThree-\thetaTwo) -- (tau3); 
	\draw[fill=white] (\thetaTwo/60,\lThree+\thetaThree-\thetaTwo) circle (2pt);
	\draw[fill=black] (tau3) circle (2pt);
	\draw[thick] (\thetaThree/60,\lFour+\thetaFour-\thetaThree) -- (tau4);	
	\draw[fill=white] (\thetaThree/60,\lFour+\thetaFour-\thetaThree) circle (2pt);
	\draw[fill=black] (tau4) circle (2pt);
	\draw[thick] (\thetaFour/60,\lOne+\thetaOne+60-\thetaFour) -- (60/60,\lOne+\thetaOne);
	\draw[fill=white] (\thetaFour/60,\lOne+\thetaOne+60-\thetaFour) circle (2pt);

	
\end{tikzpicture}
	\end{minipage}
	\caption{Alternating timetable with observed route set measure $\mathcal{M}^{RS}_{sp}(\routes_\timetable(t))$.}
	\vspace{\baselineskip}
\end{subfigure}
\begin{subfigure}[t]{\textwidth}
	\centering
	\begin{minipage}{0.35\textwidth}
		\centering
		\includegraphics{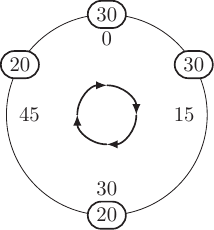}
	\end{minipage}
	\hspace{0.01\textwidth}
	\begin{minipage}{0.5\textwidth}
		\centering
		\def\lOne{30}
\def\lTwo{30}
\def\lThree{20}
\def\lFour{20}
\def\thetaOne{0}
\def\thetaTwo{10}
\def\thetaThree{30}
\def\thetaFour{50}
\begin{tikzpicture}[x=5cm,y=0.1cm]
	\draw [thick] (0,40) -- (0,10) -- (1,10);
	\foreach \x in {10, 20, 30}
	\draw[thick] (-0.5/40, \x) -- (0, \x);
	\foreach \x in {15, 20, 25, 30, 35, 40}
	\draw[thin,dotted] (0,\x) -- (1,\x);  
	\foreach \x/\xtext in {10, 20, 30, 40}
	\node [left] at (0,\x) {$\xtext$};
	\foreach \x in {0, 15, 30, 45, 60}
	\draw[thick] (\x/60, 10) -- (\x/60, 10-0.5);
	\foreach \x in {0, 5, 10, 15, 20, 25, 30, 35, 40, 45, 50, 55}
	\draw[thin,dotted] (\x/60, 40) -- (\x/60, 10);
	
	
	\node [left] at (-0.15,25) {$\mathcal{M}^{RS}_{sp}$};
	\node [below] at (0,10) {0};
	\node [below] at (0.25,10) {15};
	\node [below] at (0.5,10) {30};
	\node [below] at (0.75,10) {45};
	\node [below] at (1,10) {60};
	\node [below] at (0.5,10-4) {$t$};
	
	\coordinate (tau1) at (\thetaOne/60,\lOne);
	\coordinate (tau2) at (\thetaTwo/60,\lTwo);
	\coordinate (tau3) at (\thetaThree/60,\lThree);
	\coordinate (tau4) at (\thetaFour/60,\lFour);
	
	\draw[thick] (0,\lOne+\thetaOne-0) -- (tau1);
	\draw[fill=black] (tau1) circle (2pt);
	\draw[thick] (\thetaOne/60,\lTwo+\thetaTwo-\thetaOne) -- (tau2); 
	\draw[fill=white] (\thetaOne/60,\lTwo+\thetaTwo-\thetaOne) circle (2pt);
	\draw[fill=black] (tau2) circle (2pt);
	\draw[thick] (\thetaTwo/60,\lThree+\thetaThree-\thetaTwo) -- (tau3); 
	\draw[fill=white] (\thetaTwo/60,\lThree+\thetaThree-\thetaTwo) circle (2pt);
	\draw[fill=black] (tau3) circle (2pt);
	\draw[thick] (\thetaThree/60,\lFour+\thetaFour-\thetaThree) -- (tau4);	
	\draw[fill=white] (\thetaThree/60,\lFour+\thetaFour-\thetaThree) circle (2pt);
	\draw[fill=black] (tau4) circle (2pt);
	\draw[thick] (\thetaFour/60,\lOne+\thetaOne+60-\thetaFour) -- (60/60,\lOne+\thetaOne);
	\draw[fill=white] (\thetaFour/60,\lOne+\thetaOne+60-\thetaFour) circle (2pt);
	
	
\end{tikzpicture}
	\end{minipage}
	\caption{Non-alternating timetable with observed route set measure $\mathcal{M}^{RS}_{sp}(\routes_\timetable(t))$.}
\end{subfigure}
\caption{Optimal timetables for different route orderings.}
\label{fig:order_comp}
\end{figure}

An interesting property of optimal timetables is that \emph{for any permutation of the routes, there exists an optimal timetable with departures in that order}.
This follows from the fact that objectives~\eqref{eq:sp:opt:obj} and \eqref{eq:logit:opt:obj} are separable, and relabeling the routes does not affect the optimal solution.
Furthermore, the specific ordering of the routes plays no role in Algorithm~\ref{alg:construct}.
Figure~\ref{fig:order_comp} demonstrates this property with two optimal timetables.
The timetable at the top alternates slower and faster routes, while the timetable at the bottom follows two slow routes by two fast routes.
We obtain the managerial insight that both cyclical orderings can be optimal if the spacing between departures is optimized accordingly.

\paragraph*{Solving Problems~\eqref{formulation:sp} and \eqref{formulation:logit}}

Next, we focus on solving Problems~\eqref{formulation:sp} and \eqref{formulation:logit}.
We have already seen that the solutions can be used by Algorithm~\ref{alg:construct} to construct an optimal timetable, and the optimal objective values will be used as line plan measures.
A crucial observation is that \eqref{formulation:sp} and \eqref{formulation:logit} are convex optimization problems that can be solved efficiently.
The convexity stems from the convex timetable measures \eqref{eq:tt_sp_convex} and \eqref{eq:tt_logit_convex}, and the fact that all constraints are linear.
The $\Delta$-representation has been crucial to achieve this result, as the optimization problem for logit routing is not convex in the standard $\delta$-representation.

There are a number of tools and techniques to solve Problems~\eqref{formulation:sp} and \eqref{formulation:logit}.
Using a classical cutting plane method \citep{KelleyJr.1960-CuttingPlaneMethod}, the convex objective can be enforced through linear constraints.
This method is compatible with most mixed-integer linear programming solvers such as CPLEX and Gurobi, making it easy to integrate the new line plan measures with other problems in the same framework.
In fact, the quadratic objective of Problem~\eqref{formulation:sp} is supported directly by both solvers.
Another option is to use general non-linear optimization solvers such as Ipopt that converge to a local optimum, which is guaranteed to be globally optimal for convex problems.

When considered in isolation, the two problems can be seen as continuous non-linear resource allocation problems \citep{Patriksson2008-SurveyContinuousNonlinear}.
For shortest path routing, the resource is the period $T$ that is allocated to the gaps between departures.
For logit routing, the same resource is instead allocated to the jumps in the observed route set measure function.
\citet{Patriksson2008-SurveyContinuousNonlinear} surveys non-linear resource allocation problems, and identifies Lagrange multiplier methods as the most common solution technique.
These methods introduce a multiplier $\mu \in \mathbb{R}$ to penalize the budget constraint~\eqref{eq:sp:opt:total} or \eqref{eq:logit:opt:total} in the objective.
For a given $\mu$, each allocation $i \in \routes$ is calculated independently.
The optimal value of $\mu$ can be found with a simple line search to find the point where the independent allocations together satisfy the budget constraint.

For shortest path routing in particular, the Lagrange multiplier method reveals that every optimal solution to Problem~\eqref{formulation:sp} satisfies	$x_i(\mu) = \max\{0, \mu - l_i\}$ for some Lagrange multiplier $\mu \in \mathbb{R}$ (Appendix~\ref{app:lp_proofs}.\ref{prop:xoptstruc}).
This structure can also be observed in Figure~\ref{fig:order_comp}.
The optimal timetable only contains routes $i \in \routes$ for which $x_i > 0$, and thus $\delta_i = x_i = \mu - l_i$.
Furthermore, we have that $\tau_i = l_i$.
It follows that when the function jumps, it always jumps to $\delta_i + \tau_i = \mu$, regardless of the route.
The optimal multiplier $\mu$ is the value for which the function $\sum_{i \in \routes} x_i(\mu) = \sum_{i \in \routes} \max\{0, \mu - l_i\}$ equals $T$ (Appendix~\ref{app:lp_proofs}.\ref{col:xoptstruc:unique}).
This is a piecewise-linear function with $O(n)$ breakpoints, and bisection search can be used to find the intersection with $T$.
The corresponding complexity is $O(n \log n)$ in general, or $O(n)$ when the durations $l_i$ are already sorted \citep{Patriksson2008-SurveyContinuousNonlinear}.
Alternatively, \citet{Brucker1984-OnAlgorithmQuadratic} presents a more complicated algorithm that solves Problem~\eqref{formulation:sp} in $O(n)$ even when the routes are in arbitrary order.

\paragraph*{Line Plan Measures}

We propose to evaluate the quality of a line plan by calculating the measure of the optimal timetable.
\begin{itemize}
	\item \emph{Line Plan Measure for Shortest Path Routing or Logit Routing (general demand):}
	\begin{equation}
		\mathcal{M}^{LP}(\routes, T) = \min_{\timetableBold \in \left[0,T\right)^n}\mathcal{M}^{TT}(\routes, T, \timetableBold).
	\end{equation}
\end{itemize}
This measure inherits (strict) monotonicity from the timetable measures.
After all, when routes are reduced in duration or when new routes are added, all timetable measures improve, including the minimum.

The general measures are again specialized to uniform demand.
In this case, the measures can be calculated efficiently by solving Problem~\eqref{formulation:sp} or Problem~\eqref{formulation:logit}, respectively, as discussed in the previous section.
We obtain:
\begin{itemize}
	\item \emph{Line Plan Measure for Shortest Path Routing (uniform demand):}
	\begin{equation}
		\mathcal{M}^{LP}_{sp}(\routes, T) = \min_{\substack{\sum x_i = T \\ x\ge 0}} \frac{1}{T} \sum_{i\in \routes} \left( \frac{1}{2} x_i^2 + l_i x_i \right).
	\end{equation}
	\item \emph{Line Plan Measure for Logit Routing (uniform demand):}
	\begin{equation}
		\mathcal{M}^{LP}_{logit}(\routes, T) = \min_{\substack{\sum y_i = T \\ y\ge 0}}  \frac{1}{T} \sum_{i\in \routes} \left( \frac{1}{2} y_i^2 + l_i y_i + \frac{y_i}{\beta}\log\left(\frac{1 - e^{-\beta y_i}}{1 - e^{-\beta T}}\right) \right).
	\end{equation}
\end{itemize}


\paragraph*{Route Choice Models for Line Plans}

We again observe a natural parallel between measures and route choice models for line plans.
In the case of shortest path routing, it was observed in Section~\ref{sec:tt}, \eqref{eq:tt_routing_sp_unif}, that travelers choose route $i \in \routes$ with probability $\delta_i/T$.
In the optimal timetable, the values of $\delta$ are obtained from the $x$-variables in Problem~\eqref{formulation:sp}.
The substitution $p_i = x_i/T$ then results in the following route choice model:
\begin{itemize}
	\item \emph{Route Choice Model for Shortest Path Routing (uniform demand):}
	\begin{equation}
		\routingBold^{LP}_{sp} = \argmin_{p \in \mathcal{P}} \sum_{i \in \routes}\left( \frac{1}{2} T p_i^2 + l_i p_i\right).
	\end{equation}
\end{itemize}
Note that constraint \eqref{eq:sp:opt:total} is captured by $p \in \mathcal{P}$, which requires the probabilities to sum to one.

The route choice model expresses how travelers distribute over the routes when the timetable is optimized.
These travelers do not only care about the expected travel time $\sum_{i\in \routes} l_i p_i$, but also appreciate a small Simpson index $\sum_{i \in \routes} p^2_i$, which indicates a high diversity in route choices \citep{Simpson1949-MeasurementDiversity}.
This diversity term is given more weight when the period $T$ increases, as not diversifying the routes comes with larger consequences in terms of waiting time.
When the period decreases $T \rightarrow 0$, the route choice model simply selects the shortest path, matching the route choice model of the underlying route set.

Finally, and most surprisingly, travelers under logit routing use route $i \in \routes$ with probability $p_i = \Delta_i/T$, but only when facing the optimal timetable (Appendix~\ref{app:lp_proofs}.\ref{lemma:output_delta_logit}-\ref{prop:face_prob}).
This simple expression suggests that the $\Delta$-representation is in some sense natural for logit routing.
As the line plan indeed uses the optimal timetable, a similar substitution of $p_i = y_i/T$ can be made to obtain:
\begin{itemize}
	\item \emph{Route Choice Model for Logit Routing (uniform demand):}
	\begin{equation}
		\routingBold^{LP}_{logit} = \argmin_{p \in \mathcal{P}} \sum_{i \in \routes}\left( \frac{1}{2} T p_i^2 + l_i p_i + \frac{p_i}{\beta} \log \left(\frac{1-e^{-\beta T p_i}}{1-e^{-\beta T}}\right)\right).
	\end{equation}
\end{itemize}
Compared to shortest path routing, this route choice model includes an additional term that puts extra emphasis on avoiding probabilities that are close to zero or one.
This encourages multiple good route options to be available, regardless of when the traveler arrives.
For a decreasing period $T \rightarrow 0$, the new term tends to $\frac{1}{\beta} p_i \log(p_i)$ in the limit.
The route choice model then minimizes perceived travel time, and matches the route choice model of the underlying route set. For an increasing period $T \rightarrow \infty$, the new term approaches 0, implying that logit routing converges to shortest path routing as $T$ increases.
Finally, we note the deep connection between the line plan measures and the route choice models: they correspond to the objective value and the optimal solution to the same optimization problem.

\section{Applications}
\label{sec:ill}

Next, we demonstrate how the new framework can be used to gain practical insights and help plan real transportation systems.
Figure~\ref{fig:dutch_swiss_networks} presents portions of the Dutch and Swiss national railway networks that will be used as case studies.
Both networks operate on a periodic timetable with $T=60$ minutes, and $\beta =0.2$ is assumed for experiments with logit routing.
Route options, travel durations, and departure times are collected from the official travel planners \citep{nsreisplanner,sbbhomepage}.

\begin{figure}[th]
\centering

\tikzset{
  station/.style={ 
    draw=black, 
    fill=white, 
    inner sep=2pt, 
    font=\small,
    align=center,
    rounded corners=4pt,
  },
  route/.style={thick, black},
}

\begin{subfigure}[t]{0.49\textwidth}
\centering
\begin{tikzpicture}
\begin{axis}[
    axis lines=none,
    xtick=\empty,
    ytick=\empty,
    xmin=4.3, xmax=5.11,
    ymin=51.8, ymax=52.4,
    width=\linewidth,
    height=8cm,
    clip=false,
]

\node[station] (Rotterdam) at (axis cs:4.46944,51.92444) {Rotterdam};
\node[station] (AmsterdamC) at (axis cs:4.90000,52.37833) {Amsterdam};
\node[station] (DenHaag) at (axis cs:4.32455,52.08167) {Den Haag};
\node[station] (Utrecht) at (axis cs:5.10972,52.08917) {Utrecht};
\node[station] (Gouda) at (axis cs:4.70556,52.01750) {Gouda};
\node[station] (Schiphol) at (axis cs:4.72167,52.30089) {Schiphol};
\node[station] (Haarlem) at (axis cs:4.63889,52.38778) {Haarlem};
\node[station] (Leiden) at (axis cs:4.48222,52.16639) {Leiden};
\node[station] (Bijlmer) at (axis cs:4.94611,52.30306) {Bijlmer};

\draw[route] (AmsterdamC) -- (Schiphol);
\draw[route] (Schiphol) -- (Rotterdam);
\draw[route] (Rotterdam) -- (Gouda);
\draw[route] (Gouda) -- (Utrecht);
\draw[route] (AmsterdamC) -- (Bijlmer);
\draw[route] (Haarlem) -- (AmsterdamC);
\draw[route] (DenHaag) -- (Gouda);
\draw[route] (DenHaag) -- (Leiden);
\draw[route] (Leiden) -- (Schiphol);
\draw[route] (Bijlmer) -- (Schiphol);
\draw[route] (Utrecht) -- (Bijlmer);
\draw[route] (Leiden) -- (Utrecht);
\draw[route] (Leiden) -- (Haarlem);
\draw[route] (Schiphol) -- (Haarlem);
\draw[route] (Gouda) -- (Bijlmer);
\draw[route] (DenHaag) -- (Rotterdam);

\end{axis}
\end{tikzpicture}
\caption{Dutch network.}
\end{subfigure}
\hfill
\begin{subfigure}[t]{0.49\textwidth}
\centering
\begin{tikzpicture}
\begin{axis}[
    axis lines=none,
    xtick=\empty,
    ytick=\empty,
    xmin=6.5, xmax=8.5,
    ymin=45.5, ymax=47.5,
    width=\linewidth,
    height=8cm,
    clip=false,
]

\node[station] (Bern) at (axis cs:7.44744,46.94809) {Bern};
\node[station] (Basel) at (axis cs:7.58857,47.55960) {Basel};
\node[station] (Interlaken) at (axis cs:7.85252,46.68333) {Interlaken};
\node[station] (Lausanne) at (axis cs:6.63333,46.51667) {Lausanne};
\node[station] (Olten) at (axis cs:7.90374,47.35185) {Olten};
\node[station] (Zurich) at (axis cs:8.54169,47.37818) {Z\"urich};
\node[station] (Luzern) at (axis cs:8.30931,47.05016) {Luzern};
\node[station] (Visp) at (axis cs:7.53667,46.31861) {Visp};
\node[station] (Biel) at (axis cs:7.24336,47.13643) {Biel};

\draw[route] (Bern) -- (Olten);
\draw[route] (Bern) -- (Lausanne);
\draw[route] (Bern) -- (Biel);
\draw[route] (Olten) -- (Biel);
\draw[route] (Basel) -- (Olten);
\draw[route] (Olten) -- (Zurich);
\draw[route] (Olten) -- (Luzern);
\draw[route] (Luzern) -- (Zurich);
\draw[route] (Lausanne) -- (Visp);
\draw[route] (Lausanne) -- (Biel);
\draw[route] (Visp) -- (Bern);
\draw[route] (Visp) -- (Interlaken);
\draw[route] (Interlaken) -- (Bern);
\draw[route] (Biel) -- (Basel);
\draw[route] (Zurich) -- (Basel);
\draw[route] (Interlaken) -- (Luzern);

\end{axis}
\end{tikzpicture}
\caption{Swiss network.}
\end{subfigure}

\caption{Considered parts of the Dutch and Swiss railway networks.}
\label{fig:dutch_swiss_networks}
\end{figure}

\subsection{Line Planning}
\label{subsec:lp-app}
Line planning is one of the key strategic decisions faced by public transport operators \citep{HuismanMaroti2024-OperationsResearchNetherlands}.
Based on estimated passenger demands, the line plan determines which direct connections are available and at which frequency they are served.
To design an effective line plan, it is crucial to understand the downstream impact of the choices made at the line planning stage.
A major challenge is that the timetable is not yet known, and determining an optimal timetable at this stage is computationally prohibitive.
After all, the timetabling problem is already NP-hard by itself for just two OD pairs (see Appendix~\ref{app:lp_proofs}.\ref{prop:nphard}).
The lack of a timetable makes it difficult to predict passenger route choices, which are necessary to evaluate the true quality of service.

Previous studies have avoided the need for a timetable by making simplifying assumptions.
This include ignoring waiting time entirely, as well as assuming equidistance departures  \citep{Hartleb2022,schmidt2024planning}.
In the former case, this implies that passengers make the same decision regardless of when they arrive, while in the latter case, passengers assign the same probability to all available options.
In addition to these shortcomings, scalability remains a challenge as well.
To solve line planning problems at practical scale, Netherlands Railways relies on a genetic algorithm rather than on mixed-integer programming methods; repeatedly evaluating line plans for tens of thousands of OD pairs as part of the heuristic \citep{HuismanMaroti2024-OperationsResearchNetherlands}.

The line plan measures proposed in this paper avoid the simplifying assumptions of prior work and can directly be implemented at scale as part of existing heuristics---like the one used by Netherlands Railways.
Rather than ignoring waiting time or assuming equidistance spacing, the new measures use optimal route spacings.
Computationally, evaluating a line plan amounts to solving one (typically small) convex optimization problem for each OD pair: Problem~\eqref{formulation:sp} for shortest path routing or Problem~\eqref{formulation:logit} for logit routing.
In the case of shortest path routing, this is possible in $O(n)$ time for $n$ route options (see Section~\ref{sec:lp}).
The problems are solved for each OD pair independently and can be solved in parallel if needed, leading to a scalable solution with a consistent theoretical foundation.

An important question that remains is whether the line plan measures accurately capture the properties of the timetable that will eventually be implemented downstream.
After all, each OD pair is optimized independentenly and interactions between OD pairs are ignored at the line planning stage.
To empirically test the quality of this approximation, we compare two values for the Dutch and Swiss rail networks: the line plan measure that is calculated based only on route options and corresponding travel times, and the timetable measure of the \emph{actual} timetable that was implemented in practice.
If these values are close, this indicates that the line plan measure provides a good approximation of the service quality of the system that will eventually be implemented, despite not yet having access to the timetable.

\begin{figure}[ht]
\centering

\newcommand{\axiswidth}{7.5cm}
\newcommand{\axisheight}{7.5cm}

\begin{subfigure}[t]{0.49\textwidth}
\centering
\begin{tikzpicture}
\begin{axis}[
    width=\axiswidth,
    height=\axisheight,
    xlabel={$\mathcal{M}^{LP}_{sp}$},
    ylabel={$\mathcal{M}^{TT}_{sp}$},
    every axis plot/.append style={thick},
    tick label style={font=\small},
    label style={font=\normalsize},
    legend style={at={(0.05,0.95)},anchor=north west, font=\small},
    legend cell align={left},
    xmin=5, xmax=75, ymin=5, ymax=75
]
\addplot[only marks, mark=*, mark size=1.5pt] table [x=m_LP_SP, y=m_TT_SP, col sep=comma] {scatter.csv};
\addplot[red, very thin] table [x=x, y=y, col sep=comma] {line.csv};
\end{axis}
\end{tikzpicture}
\caption{Dutch Network: Shortest Path}
\end{subfigure}
\hfill
\begin{subfigure}[t]{0.49\textwidth}
\centering
\begin{tikzpicture}
\begin{axis}[
    width=\axiswidth,
    height=\axisheight,
    xlabel={$\mathcal{M}^{LP}_{logit}$},
    ylabel={$\mathcal{M}^{TT}_{logit}$},
    every axis plot/.append style={thick},
    tick label style={font=\small},
    label style={font=\normalsize},
    legend style={at={(0.05,0.95)},anchor=north west, font=\small},
    legend cell align={left},
    xmin=5, xmax=75, ymin=5, ymax=75
]
\addplot[only marks, mark=*, mark size=1.5pt] table [x=m_LP_Logit, y=m_TT_Logit, col sep=comma] {scatter.csv};
\addplot[red, very thin] table [x=x, y=y, col sep=comma] {line.csv};
\end{axis}
\end{tikzpicture}
\caption{Dutch Network: Logit ($\beta=0.2$)}
\end{subfigure}

\vspace{0.5cm}

\begin{subfigure}[t]{0.49\textwidth}
\centering
\begin{tikzpicture}
\begin{axis}[
    width=\axiswidth,
    height=\axisheight,
    xlabel={$\mathcal{M}^{LP}_{sp}$},
    ylabel={$\mathcal{M}^{TT}_{sp}$},
    every axis plot/.append style={thick},
    tick label style={font=\small},
    label style={font=\normalsize},
    legend style={at={(0.05,0.95)},anchor=north west, font=\small},
    legend cell align={left},
    xmin=30, xmax=180, ymin=30, ymax=180
]
\addplot[only marks, mark=*, mark size=1.5pt] table [x=m_LP_SP, y=m_TT_SP, col sep=comma] {sbb.csv};
\addplot[red, very thin] table [x=x, y=y, col sep=comma] {line.csv};
\end{axis}
\end{tikzpicture}
\caption{Swiss Network: Shortest Path}
\end{subfigure}
\hfill
\begin{subfigure}[t]{0.49\textwidth}
\centering
\begin{tikzpicture}
\begin{axis}[
    width=\axiswidth,
    height=\axisheight,
    xlabel={$\mathcal{M}^{LP}_{logit}$},
    ylabel={$\mathcal{M}^{TT}_{logit}$},
    every axis plot/.append style={thick},
    tick label style={font=\small},
    label style={font=\normalsize},
    legend style={at={(0.05,0.95)},anchor=north west, font=\small},
    legend cell align={left},
    xmin=30, xmax=180, ymin=30, ymax=180
]
\addplot[only marks, mark=*, mark size=1.5pt] table [x=m_LP_Logit, y=m_TT_Logit, col sep=comma] {sbb.csv};
\addplot[red, very thin] table [x=x, y=y, col sep=comma] {line.csv};
\end{axis}
\end{tikzpicture}
\caption{Swiss Network: Logit ($\beta=0.2$)}
\end{subfigure}

\caption{Comparison of line plan and timetable measures for Dutch and Swiss rail networks.}
\label{fig:scatter_plots}
\end{figure}

Figure~\ref{fig:scenarios} demonstrates that the line plan measures predict the quality of the real systems almost perfectly.
For all station pairs, the figure compares the line plan measure (without access to the timetable) and the timetable measure (using the real timetable).
Across both countries and both routing approaches, the measures differ only by about 1\% on average, indicating that the line plan measure is about 99\% accurate in predicting the measure of the eventual timetable.
This close match can likely be attributed to the fact that the line plan measures anticipate that timetables will be optimized downstream, and both public transport operators indeed produced timetables that are close to optimal.
It also suggests that the new line plan measures have major practical value in capturing the impact of timetables at the line planning stage, something that is otherwise very hard to do.

Another advantage is that the new measures can provide realistic insights into the expected passenger flows, already at the line planning stage.
After all, the distribution of the passengers over the routes is obtained as a simple byproduct when calculating the measure (see Section~\ref{sec:lp}).
To evaluate the accuracy of these predictions for the Dutch and the Swiss rail networks, we calculate the total variation distance $\frac{1}{2} \lVert p^{LP} - p^{TT} \rVert_1$ between the routing predicted by the line plan measure ($p^{LP}$) and the routing that is obtained when the real timetable is available ($p^{TT}$).
While the individual passenger flows are harder to estimate than the total service quality, the predictions are still remarkably accurate: the total variation distance is under 9\% on average for both the Dutch and the Swiss network, meaning that 91\% of flows are predicted correctly.
Accurate flow predictions in turn open the door for capacity planning and congestion analysis at the line planning stage.

In addition to integrating the new line plan measures into existing heuristics to obtain the benefits discussed above, it would also be very interesting to integrate them into existing mixed-integer programming models.
The line plan measures avoid the simplifying assumptions that are currently used and even open up new opportunities from a computational perspective: the convex structure of Problems~\eqref{formulation:sp} and~\eqref{formulation:logit} naturally lends itself to decomposition schemes such as  Benders decomposition, offering a potential path forward to solving these integrated models at scale.

\subsection{Analyzing Timetables}

The new measures can also be used to analyze timetables and to point out potential inefficiencies.
By design, the line plan measure equals the timetable measure of the optimal timetable for that connection.
Therefore, differences between the two measures may indicate a suboptimal timetable in parts of the network.

As an example, Figure~\ref{fig:actandopt} highlights the most significant outlier from Figure~\ref{fig:scatter_plots}a: the connection between Rotterdam and Bijlmer. Here, the line plan measure is 51.1 minutes, while the timetable measure is 54.6 minutes. This OD pair has three fast routes (durations of 39, 40, and 48 minutes) and three slower routes (61, 61, and 63 minutes). In the optimal timetable, the fast routes are evenly spaced. However, in the actual timetable, the fast routes are clustered within a 30-minute window. As a result, travelers have longer waiting times, or resort to slower routes, increasing average travel time which causes the discrepancy.
\begin{figure}[t]
\def\lOne{63}
\def\lTwo{48}
\def\lThree{39}
\def\lFour{40}
\def\lFive{61}
\def\lSix{61}
\def\thetaOne{5}
\def\thetaTwo{22}
\def\thetaThree{41}
\def\thetaFour{51}
\def\thetaFive{54}
\def\thetaSix{24.5}

\centering
\begin{subfigure}[b]{0.35\textwidth}
\centering
\begin{tikzpicture}[>=latex]
	\coordinate (O) at (0,0);
	\def\radius{1.7cm}
	\def\smallRadius{0.5cm}
	  
	\coordinate (left) at (180:\smallRadius);
	\coordinate (right) at (0:\smallRadius); 
	\coordinate (top) at (90:\smallRadius);
	\coordinate (bottom) at (-90:\smallRadius);
	\draw[->,thick,out=90,in=180] (left) to (top);
	\draw[->,thick,out=0,in=90] (top) to (right);
	\draw[->,thick,out=-90,in=0] (right) to (bottom);
	\draw[->,thick,out=180,in=-90] (bottom) to (left);
	
	\draw (O) circle[radius=\radius];
	
	\coordinate (0t) at (90:\radius); 
	\fill (0t)  (90: \radius-0.32cm) node{$0$};
	\coordinate (180t) at (-90:\radius); 
	\fill (180t)  (-90: \radius-0.35cm) node{30};
	\coordinate (15t) at (0:\radius); 
	\fill (15t)  (0: \radius-0.39cm) node{15};
	\coordinate (45t) at (180:\radius); 
	\fill (45t)  (180: \radius-0.39cm) node{45};
	
	\coordinate (e1) at (90-\thetaOne*6:\radius);
	\coordinate (e2) at (90-\thetaTwo*6:\radius);
	\coordinate (e3) at (90-\thetaThree*6:\radius);
	\coordinate (e4) at (90-\thetaFour*6:\radius);
	\coordinate (e5) at (90-\thetaFive*6:\radius);
	\coordinate (e6) at (90-\thetaSix*6:\radius);

    \node[vertex, fill=white] at (e6) {\lSix};
	\node[vertex, fill=white] at (e1) {\lOne};
	\node[vertex, fill=white] at (e2) {\lTwo};
	\node[vertex, fill=white] at (e3) {\lThree};
	\node[vertex, fill=white] at (e4) {\lFour};
	\node[vertex, fill=white] at (e5) {\lFive};

	\node at (0,-2.3) {\textcolor{white}{text}};
\end{tikzpicture}
\caption{Actual.}
\label{fig:actualRtdAsb}
\end{subfigure}
\hspace{0.01\textwidth}
\def\lOne{61}
\def\thetaOne{24}
\def\thetaTwo{23}
\def\thetaThree{47}
\def\thetaFour{8}
\def\thetaFive{9}
\begin{subfigure}[b]{0.50\textwidth}
\centering
\begin{tikzpicture}[>=latex]
	\coordinate (O) at (0,0);
	\def\radius{1.7cm}
	\def\smallRadius{0.5cm}
	  
	\coordinate (left) at (180:\smallRadius);
	\coordinate (right) at (0:\smallRadius); 
	\coordinate (top) at (90:\smallRadius);
	\coordinate (bottom) at (-90:\smallRadius);
	\draw[->,thick,out=90,in=180] (left) to (top);
	\draw[->,thick,out=0,in=90] (top) to (right);
	\draw[->,thick,out=-90,in=0] (right) to (bottom);
	\draw[->,thick,out=180,in=-90] (bottom) to (left);
	
	\draw (O) circle[radius=\radius];
	
	\coordinate (0t) at (90:\radius); 
	\fill (0t)  (90: \radius-0.32cm) node{$0$};
	\coordinate (180t) at (-90:\radius); 
	\fill (180t)  (-90: \radius-0.35cm) node{30};
	\coordinate (15t) at (0:\radius); 
	\fill (15t)  (0: \radius-0.39cm) node{15};
	\coordinate (45t) at (180:\radius); 
	\fill (45t)  (180: \radius-0.39cm) node{45};
	
	\coordinate (e1) at (90-\thetaOne*6:\radius);
	\coordinate (e2) at (90-\thetaTwo*6:\radius);
	\coordinate (e3) at (90-\thetaThree*6:\radius);
    \coordinate (e5) at (90-\thetaFive*6:\radius);
	\coordinate (e4) at (90-\thetaFour*6:\radius);

	\node[vertex, fill=white] at (e1) {\lOne};
	\node[vertex, fill=white] at (e2) {\lTwo};
	\node[vertex, fill=white] at (e3) {\lThree};
    \node[vertex, fill=white] at (e5) {\lFive};
	\node[vertex, fill=white] at (e4) {\lFour};

	\node at (0,-2.3) {\textcolor{white}{text}};
\end{tikzpicture}
\caption{Optimal.}
\label{fig:optRtdAsb}
\end{subfigure}

\caption{Actual and optimal timetable for Rotterdam - Bijlmer (using shortest path).}
\label{fig:actandopt}
\end{figure}

While individual inefficiencies cannot always be resolved in the context of an integrated timetable, the measures proposed in this paper provide an easy way to identify potential weaknesses.
This is especially important when timetabling is partially manual and thousands of OD pairs are involved.
In fact, \cite{HuismanMaroti2024-OperationsResearchNetherlands} identify the timetabling process at Netherlands Railways as a key focus area that needs more support in the planning pipeline.
The new measures also enable analysis from a fairness perspective, e.g., by mapping whether the timetable in different parts of the country is equally close to optimal.

\subsection{Managerial Insights}

Finally, we use the new framework to obtain managerial insights about 1) the importance of using the right model at the right time, 2) when shortest path routing is all you need, and 3) the consequences of using inconsistent models.
To this end we zoom in on the the connection between Tilburg and Eindhoven in the Netherlands and we consider an example to demonstrate how inconsistency may impact decision making.

\paragraph{Using the Right Model at the Right Time}

To demonstrate the importance of using the right model at the right time, we examine in detail the measures and route choice models on the rail connection between the Dutch cities of Tilburg (Tb) and Eindhoven (Ehv). In the 2025 timetable, Netherlands Railways operates one faster Intercity (IC) train, and one slower Sprinter (Spr) train between these cities every 30 minutes, with approximately equidistant departures. The corresponding route set, timetable, and line plan  are visualized in Figure~\ref{fig:concepts_I}.

Starting from this baseline, the following figures show how the routings and measures change when varying the route duration of the Sprinter train $l_\text{Spr}$ (Figure~\ref{fig:illustration1}), the period $T$ (Figure~\ref{fig:illustration2}), and the departure time of the Sprinter train $\timetable_\text{Spr}$ (Figure~\ref{fig:illustration3}).

\begin{figure}[p]
	\centering
	\begin{subfigure}[t]{0.165\textwidth}
		\centering
		\fbox{\includegraphics[height=2.2cm]{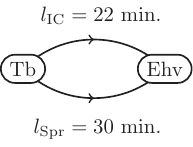}}
		\caption{Route Set.}
	\end{subfigure}
	\hspace{0.03\textwidth}
	\begin{subfigure}[t]{0.35\textwidth}
		\centering
		\fbox{\includegraphics[height=2.2cm]{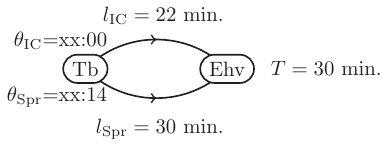}}
		\caption{Timetable.}
	\end{subfigure}
	\hspace{0.03\textwidth}
	\begin{subfigure}[t]{0.35\textwidth}
		\centering
		\fbox{\includegraphics[height=2.2cm]{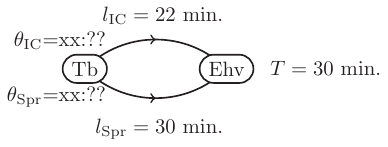}}
		\caption{Line Plan.}
	\end{subfigure}
	\caption{Baseline route set, timetable, and line plan for sensitivity analysis.}
	\label{fig:concepts_I}
\end{figure}

\begin{figure}[p]
\centering
    \begin{subfigure}[b]{0.4\textwidth}
    \begin{tikzpicture}[scale=0.68]
\begin{axis}[cycle list name=black white,xlabel={Duration Sprinter},
  ylabel={Probability Sprinter},
xmin = 0, xmax = 40,legend style={at={(1.02,1)},anchor=north west}]
 \addplot+[mark=none,domain = 0:22] {1};
             \draw[fill=white] (22,1) circle (2pt);
            \draw[fill=black] (22,0) circle (2pt);
\addplot+[mark=none, dashed]  table [x=Traveltime2, y=pLogit_RS, col sep=comma] {data1.csv};
\addplot+[mark=none, dotted]  table [x=Traveltime2, y=pLogit_RS2, col sep=comma] {data1.csv};
\addplot+[mark=none,domain = 22:40] {0};
\end{axis}
\end{tikzpicture}
    \caption{Route Set Routing.}
    \label{fig:illustration1a}
    \end{subfigure}
	\vspace{\baselineskip}
        \begin{subfigure}[b]{0.45\textwidth}
    \begin{tikzpicture}[scale=0.68]
\begin{axis}[cycle list name=black white,xlabel={Duration Sprinter},
  ylabel={$\mathcal{M}^{RS}$ (min.)},
xmin = 0, xmax = 40,legend style={at={(1.02,1)},anchor=north west},ymin=0,ymax=40]
\addplot+[mark=none]  table [x=Traveltime2, y=mSP_RS, col sep=comma] {data1.csv};
\addplot+[mark=none, dashed]  table [x=Traveltime2, y=mLogit_RS, col sep=comma] {data1.csv};
\addplot+[mark=none, dotted]  table [x=Traveltime2, y=mLogit_RS2, col sep=comma] {data1.csv};
\legend{Shortest Path,Logit ($\beta=0.2$),Logit ($\beta=0.1$)}
\end{axis}
\end{tikzpicture}
    \caption{Route Set Measure.}
    \label{fig:illustration1b}
    \end{subfigure}
	\vspace{\baselineskip}
    \begin{subfigure}[b]{0.4\textwidth}
    \begin{tikzpicture}[scale=0.68]
\begin{axis}[cycle list name=black white,xlabel={Duration Sprinter},
  ylabel={Probability Sprinter},
xmin = 0, xmax = 40]
 \addplot+[mark=none,domain = 0:7] {1};
    \draw[fill=white] (7,1) circle (2pt);
            \draw[fill=black] (8,0.467) circle (2pt);
            \draw[fill=black] (38,0) circle (2pt);
\addplot+[mark=none, dashed]  table [x=Traveltime2, y=pLogit_TT, col sep=comma] {data1.csv};
\addplot+[mark=none, dotted]  table [x=Traveltime2, y=pLogit_TT2, col sep=comma] {data1.csv};
\addplot+[mark=none,domain = 8:38] {0.467};
\addplot+[mark=none,domain = 38:40] {0};
    \draw[fill=white] (38,0.467) circle (2pt);
\end{axis}
\end{tikzpicture}
    \caption{Timetable Routing.}
    \label{fig:illustration1c}
    \end{subfigure}
        \begin{subfigure}[b]{0.45\textwidth}
    \begin{tikzpicture}[scale=0.68]
\begin{axis}[cycle list name=black white,xlabel={Duration Sprinter},
  ylabel={$\mathcal{M}^{TT}$ (min.)},
xmin = 0, xmax = 40,legend style={at={(0.97,0.03)},anchor=south east},ymin=0,ymax=40]
\addplot+[mark=none]  table [x=Traveltime2, y=mSP_TT, col sep=comma] {data1.csv};
\addplot+[mark=none, dashed]  table [x=Traveltime2, y=mLogit_TT, col sep=comma] {data1.csv};
\addplot+[mark=none, dotted]  table [x=Traveltime2, y=mLogit_TT2, col sep=comma] {data1.csv};
\end{axis}
\end{tikzpicture}
    \caption{Timetable Measure.}
    \label{fig:illustration1d}
    \end{subfigure}
     \begin{subfigure}[b]{0.4\textwidth}
    \begin{tikzpicture}[scale=0.68]
\begin{axis}[cycle list name=black white,xlabel={Duration Sprinter},
  ylabel={Probability Sprinter},
xmin = 0, xmax = 40,ymin = -0.1, ymax =1.1]
\addplot+[mark=none]  table [x=Traveltime2, y=pSP_LP, col sep=comma] {data1.csv};
\addplot+[mark=none, dashed]  table [x=Traveltime2, y=pLogit_LP, col sep=comma] {data1.csv};
\addplot+[mark=none, dotted]  table [x=Traveltime2, y=pLogit_LP2, col sep=comma] {data1.csv};
\end{axis}
\end{tikzpicture}
    \caption{Line Plan Routing.}
    \label{fig:illustration1e}
    \end{subfigure}
        \begin{subfigure}[b]{0.45\textwidth}
    \begin{tikzpicture}[scale=0.68]
\begin{axis}[cycle list name=black white,xlabel={Duration Sprinter},
  ylabel={$\mathcal{M}^{LP}$ (min.)},
xmin = 0, xmax = 40,legend style={at={(0.97,0.03)},anchor=south east}, ymin=0,ymax=40]
\addplot+[mark=none]  table [x=Traveltime2, y=mSP_LP, col sep=comma] {data1.csv};
\addplot+[mark=none, dashed]  table [x=Traveltime2, y=mLogit_LP, col sep=comma] {data1.csv};
\addplot+[mark=none, dotted]  table [x=Traveltime2, y=mLogit_LP2, col sep=comma] {data1.csv};
\end{axis}
\end{tikzpicture}
    \caption{Line Plan Measure.}
	\label{fig:illustration1f}
    \end{subfigure}
    \caption{
    	Sensitivity analysis for varying route duration $l_\text{Spr}$ in Figure~\ref{fig:concepts_I}.
	}
    \label{fig:illustration1}
\end{figure}

{
	\addtocounter{figure}{-1}%
	\renewcommand{\thefigure}{\ref{fig:concepts_I}}
	\begin{figure}[p]
		\centering
		\begin{subfigure}[t]{0.165\textwidth}
			\centering
			\fbox{\includegraphics[height=2.2cm]{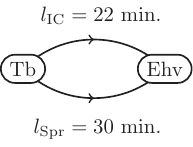}}
			\caption{Route Set.}
		\end{subfigure}
		\hspace{0.03\textwidth}
		\begin{subfigure}[t]{0.35\textwidth}
			\centering
			\fbox{\includegraphics[height=2.2cm]{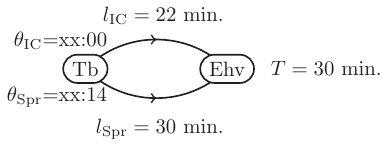}}
			\caption{Timetable.}
		\end{subfigure}
		\hspace{0.03\textwidth}
		\begin{subfigure}[t]{0.35\textwidth}
			\centering
			\fbox{\includegraphics[height=2.2cm]{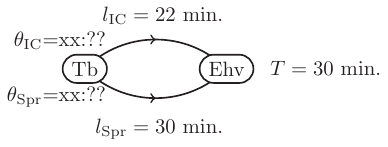}}
			\caption{Line Plan.}
		\end{subfigure}
		\caption{Baseline route set, timetable, and line plan for sensitivity analysis. (repeated)}
	\end{figure}
}

\begin{figure}[p]
\centering
    \begin{subfigure}[b]{0.4\textwidth}
    \begin{tikzpicture}[scale=0.68]
\begin{axis}[cycle list name=black white,xlabel={Period $T$ (min.)},
  ylabel={Probability Sprinter},
xmin = 0, xmax = 60,legend style={at={(1.02,1)},anchor=north west}]
 \addplot+[mark=none,domain = 0:16] {0};
   
\addplot+[mark=none, dashed]  table [x=T, y=pLogit_TT, col sep=comma] {data2.csv};
\addplot+[mark=none, dotted]  table [x=T, y=pLogit_TT2, col sep=comma] {data2.csv};
\addplot+[mark=none,domain = 16:60] {0.467};
\draw[fill=black] (16,0) circle (2pt);
 \draw[fill=white] (16,0.467) circle (2pt);
\end{axis}
\end{tikzpicture}
    \caption{Timetable Routing.}
    \label{fig:illustration2a}
    \end{subfigure}
	\vspace{\baselineskip}
        \begin{subfigure}[b]{0.45\textwidth}
    \begin{tikzpicture}[scale=0.68]
\begin{axis}[cycle list name=black white,xlabel={Period $T$ (min.)},
  ylabel={$\mathcal{M}^{TT}$ (min.)},
xmin = 0, xmax = 60,legend style={at={(1.02,1)},anchor=north west}]
\addplot+[mark=none]  table [x=T, y=mSP_TT, col sep=comma] {data2.csv};
\addplot+[mark=none, dashed]  table [x=T, y=mLogit_TT, col sep=comma] {data2.csv};
\addplot+[mark=none, dotted]  table [x=T, y=mLogit_TT2, col sep=comma] {data2.csv};
\legend{Shortest Path,Logit ($\beta=0.2$),Logit ($\beta=0.1$)}
\end{axis}
\end{tikzpicture}
    \caption{Timetable Measure.}
    \end{subfigure}
     \begin{subfigure}[b]{0.4\textwidth}
    \begin{tikzpicture}[scale=0.68]
\begin{axis}[cycle list name=black white,xlabel={Period $T$ (min.)},
  ylabel={Probability Sprinter},
xmin = 0, xmax = 60, ymax =0.55,legend style={at={(0.97,0.03)},anchor=south east}]
\addplot+[mark=none]  table [x=T, y=pSP_LP, col sep=comma] {data2.csv};
\addplot+[mark=none, dashed]  table [x=T, y=pLogit_LP, col sep=comma] {data2.csv};
\addplot+[mark=none, dotted]  table [x=T, y=pLogit_LP2, col sep=comma] {data2.csv};
\end{axis}
\end{tikzpicture}
    \caption{Line Plan Routing.}
    \end{subfigure}
        \begin{subfigure}[b]{0.45\textwidth}
    \begin{tikzpicture}[scale=0.68]
\begin{axis}[cycle list name=black white,xlabel={Period $T$ (min.)},
  ylabel={$\mathcal{M}^{LP}$ (min.)},
xmin = 0, xmax = 60,legend style={at={(0.97,0.03)},anchor=south east}]
\addplot+[mark=none]  table [x=T, y=mSP_LP, col sep=comma] {data2.csv};
\addplot+[mark=none, dashed]  table [x=T, y=mLogit_LP, col sep=comma] {data2.csv};
\addplot+[mark=none, dotted]  table [x=T, y=mLogit_LP2, col sep=comma] {data2.csv};
\end{axis}
\end{tikzpicture}
    \caption{Line Plan Measure.}
    \end{subfigure}
    \caption{Sensitivity analysis for varying period $T$ in Figure~\ref{fig:concepts_I}. The fixed timetable maintains approximately equidistant departures $\theta_\text{IC}=0$, $\theta_\text{Spr}=(14/30)T$.}
    \label{fig:illustration2}
\end{figure}

\begin{figure}[p]
\centering
    \begin{subfigure}[b]{0.4\textwidth}
    \begin{tikzpicture}[scale=0.68]
\begin{axis}[cycle list name=black white,xlabel={Departure Time Sprinter},
  ylabel={Probability Sprinter},
xmin = 0, xmax = 30,legend style={at={(1.02,1)},anchor=north west}]
\addplot+[mark=none, domain=0:21] {x/30};
\draw[fill=white] (21,21/30) circle (2pt);
\addplot+[mark=none, dashed]  table [x=q2, y=pLogit_TT, col sep=comma] {data3.csv};
\addplot+[mark=none, dotted]  table [x=q2, y=pLogit_TT2, col sep=comma] {data3.csv};
\addplot+[mark=none,domain = 21:30] {0};
\draw[fill=black] (21,0) circle (2pt);
\end{axis}
\end{tikzpicture}
    \caption{Timetable Routing.}
    \label{fig:illustration3a}
    \end{subfigure}
        \begin{subfigure}[b]{0.45\textwidth}
    \centering
    \begin{tikzpicture}[scale=0.68]
    \centering
\begin{axis}[cycle list name=black white,xlabel={Departure Time Sprinter},
  ylabel={$\mathcal{M}^{TT}$ (min.)},
xmin = 0, xmax = 30,legend style={at={(1.02,1)},anchor=north west}]
\addplot+[mark=none]  table [x=q2, y=mSP_TT, col sep=comma] {data3.csv};
\addplot+[mark=none, dashed]  table [x=q2, y=mLogit_TT, col sep=comma] {data3.csv};
\addplot+[mark=none, dotted]  table [x=q2, y=mLogit_TT2, col sep=comma] {data3.csv};
\legend{Shortest Path,Logit ($\beta=0.2$),Logit ($\beta=0.1$)}
\end{axis}
\end{tikzpicture}
    \caption{Timetable Measure.}
    \label{fig:illustration3b}
    \end{subfigure}
    \caption{Sensitivity analysis for varying departure time $\theta_\text{Spr}$ in Figure~\ref{fig:concepts_I}.}
    \label{fig:illustration3}
\end{figure}

The figures reveal that route choice and service quality can differ considerably between route sets, timetables, and line plans, demonstrating that it is important to use the right model at the right stage of planning.
In Figure~\ref{fig:illustration1}, for example, travelers under shortest path routing who are faced with a static route set will all choose the Sprinter train if the duration is less than 22 minutes.
In the timetable setting, this only happens when the Sprinter is so fast that it is the best option regardless of the arrival time (duration less than 8 minutes).
In general, Figure~\ref{fig:illustration1} shows clear differences between route sets (which ignore periodicity), and timetables and line plans (which include periodicity), both in terms of route choice and service quality.

The difference between timetable measures and line plan measures is less obvious from Figure~\ref{fig:illustration1}, but this strongly depends on the quality of the timetable.
The actual timetable is relatively close to optimal, and therefore the timetable measure is relatively close to the line plan measure (recall that the line plan measure is defined as the timetable measure of the optimal timetable).
The situation is different in Figure~\ref{fig:illustration3b}, which changes the timetable by varying the departure time of the Sprinter train.
Each curve shows the timetable measure for different timetables, and the minimum value is equal to the line plan measure (e.g., $32.967$ at $\theta_2=11$ for shortest path).
As the timetable moves further away from optimal, the difference between the timetable measure and the line plan measure becomes apparent, again emphasizing that it is important to use routings and measures tailored to the form of public transport that one considers.

\paragraph{When Shortest Path is All You Need}

In line with the literature on travel demand analysis, the smooth curves provided by the logit model suggest that it is a much more reasonable model for travel behavior than the shortest path model. Route choice based on shortest path yields large jumps in routing probabilities in Figures~\ref{fig:illustration1a}, \ref{fig:illustration1c}, \ref{fig:illustration2a} and \ref{fig:illustration3a}, which translate to kinks in the corresponding service quality measures.
In practice, routings and service quality are expected to vary smoothly when the system's parameters change, which better matches the logit model.

Another strength of the logit model is that it provides better guidance for the local improvement of existing timetables.
Consider the timetables in Figure~\ref{fig:illustration3}, for example.
For $\theta_\text{Spr} > 22$ minutes, the timetable is such that the Sprinter train is never the shortest path.
This results in a flat region of the shortest path curve in Figure~\ref{fig:illustration3b}, which provides no guidance on how to improve the timetable.
For logit the curve is smooth, and the gradient suggests moving the Sprinter up in the schedule, towards the optimal departure time.

While logit is the model of choice to accurately capture traveler behavior, the benefit of logit over shortest path seems to diminish as we move from route sets to timetables to line plans.
This can be seen in Figure~\ref{fig:illustration1}, for example, where the curves for shortest path and logit get closer to each other as we move down the plots, both for routings and for measures. We showed in Section~\ref{sec:lp} that shortest and logit yield the same results for line plans in the limit as $T\rightarrow \infty$, but this behavior already emerges for relatively small $T$, and also for timetables.

We conjecture that is due to the following two effects.
First, when there is a clear shortest path, the shortest path model and the logit model give similar results.
This can be seen in Figure~\ref{fig:illustration1a} for example, where all models assign a high probability to the Sprinter train when it is fast and a low probability when it is slow, but the models disagree when there is no clear preference around $l_\text{Spr}=22$.
Second, when travel demand is uniformly distributed, it is likely that many travelers have a clearly preferred route based on their preferred departure time, especially if the departures are spaced out well.
As such, the \emph{average} traveler is less affected by the difference between shortest path and logit.
Timetable routings and measures are defined based on the average traveler, which explains why the curves get closer together.
Line plans further strengthen this effect by optimizing the timetable to intentionally give travelers a clearly preferred choice, resulting in the smallest difference between shortest path and logit, as seen in Figures~\ref{fig:illustration1e} and \ref{fig:illustration1f}.

To validate empirically whether the difference between shortest path and logit diminishes for line plans, we once again consider the Dutch and Swiss networks introduced at the beginning of Section~\ref{sec:ill}.
The total variation distance $\frac{1}{2} \lVert p_{sp} - p_{logit} \rVert_1$ is calculated between the shortest path and logit routings to quantify how the results differ.
Figure~\ref{fig:boxplot} shows boxplots of this distance across all OD pairs, for both the Dutch and the Swiss network. It is clear that the differences between routing models are substantially smaller for line plans than for timetables, suggesting that the benefits of using a logit model indeed diminish for line planning.

\begin{figure}
\centering
\begin{tikzpicture}
\begin{axis}[
    boxplot/draw direction=y,
    ylabel={$\frac{1}{2}\left\lVert \routingBold_{sp} - \routingBold_{logit}\right\rVert_1$},
    xtick={1,2, 3, 4},
    xticklabels={{TT, Dutch}, {LP, Dutch}, {TT, Swiss}, {LP, Swiss}},
    ymajorgrids=true,
    width=14cm,
    height=7cm,
    cycle list/Set1-4
]

\addplot+[
    boxplot,
] table[y index=6,col sep=comma] {boxplotNS.csv};

\addplot+[
    boxplot,
] table[y index=7,col sep=comma] {boxplotNS.csv};

\addplot+[
    boxplot,
] table[y index=6,col sep=comma] {boxplotSBB.csv};

\addplot+[
    boxplot,
] table[y index=7,col sep=comma] {boxplotSBB.csv};

\end{axis}
\end{tikzpicture}
\caption{Comparison of the total variation distance between shortest path routing and logit routing for both timetables (TT) and line plans (LP).}
\label{fig:boxplot}
\end{figure}

The observations above suggest that \emph{it is not always necessary for public transport planners to accurately model travel behavior}.
Especially for high-level planning such as line planning, the shortest path model may be all you need to make good decisions.
This is of major managerial importance, as the shortest path model is easier to integrate in existing systems, and does not require an estimation of the $\beta$ parameter.
The framework in this paper allows practitioners to get a sense of how much shortest path and logit measures differ on their network, before investing in more complicated optimization methods that may not be necessary.

\paragraph{Negative Consequences of Inconsistent Models}

Even if there is no significant difference between the shortest path and logit models for line planning, \emph{combining the models in an inconsistent way can have significant negative consequences.}
This behavior will be demonstrated with an example below.
It was already shown in the introduction that inconsistent models can lead to counter-intuitive results for route sets, but it is surprising to see that this behavior extends all the way to line plans, which are much less affected by the choice of model.
The main purpose of this section is to serve as a warning that \emph{being consistent may be much more important than being correct}.

\begin{figure}[t]
	\centering
	\begin{subfigure}[b]{0.4\textwidth}
		\centering\includegraphics{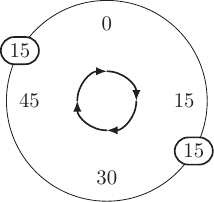}
		\caption{Original (two routes).}
	\end{subfigure}
	\vspace{\baselineskip}
	\begin{subfigure}[b]{0.45\textwidth}
		\centering
		\includegraphics{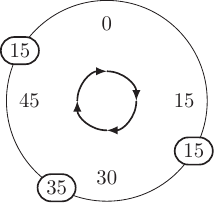}
		
		\caption{Option 1: In-between.}
	\end{subfigure}
	\begin{subfigure}[b]{0.4\textwidth}
		\centering\includegraphics{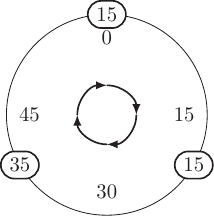}
		
		\caption{Option 2: Equidistant.}
	\end{subfigure}
	\begin{subfigure}[b]{0.45\textwidth}
		\centering\includegraphics{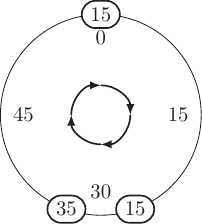}
		\caption{Option 3: Optimal.}
	\end{subfigure}
	\caption{Original timetable and three options to include a new route.}\label{fig:scenarios}
\end{figure}

The example starts from an original timetable that is depicted in Figure~\ref{fig:scenarios}a.
The timetable uses a period of one hour, and features two routes with a duration of 15 minutes that are spread out evenly.
Now suppose that the opportunity arises to add an additional route with a duration of 35 minutes.
There are a number of reasonable ways to do so:
\begin{itemize}
	\item Option 1: Place the new route directly in-between the other departures (Figure~\ref{fig:scenarios}b);
	\item Option 2: Use an equidistant timetable with three departures (Figure~\ref{fig:scenarios}c);
	\item Option 3: Use the optimal timetable (Figure~\ref{fig:scenarios}d).
\end{itemize}
The optimal timetable is calculated for the shortest path model with the methods developed in this paper.
As expected, the logit-optimal timetable is almost identical to the shortest path-optimal timetable, to the extent that there is no benefit in treating these options separately.

Table~\ref{tbl:measuresScen} summarizes how three line planners evaluate the different options.
The first planner uses the shortest path measure, the second planner uses the logit measure with $\beta=0.1$, and the third planner uses an inconsistent measure that combines logit routing with travel time evaluation (instead of the consistent \emph{perceived} travel time evaluation).

When the shortest path measure is used, adding in a slow route does not necessarily improve the service quality.
Option 1 does not affect the measure, as travelers can simply ignore the new route and wait for the faster routes they prefer.
Option 2 worsens the measure, as it ruins the spacing of the fast routes without providing any benefit.
However, the measure indicates that a strict benefit can be obtained when the routes are spaced out appropriately, as in Option 3.
In the logit case, all options are better than the original timetable.
This is due to the very low value of $\beta$, which assumes that travelers really care about having more options.
But even if this assumption were incorrect, the relative ranking between Options 1-3 is the same as for shortest path.
Planners are still encouraged to implement Option 1 over Option 2, which would not hurt in the shortest path case, or add the new route and optimize the timetable, which would result in a similar timetable and better service quality in both cases.

This is in stark contrast with the inconsistent measure: \emph{adding the new route always makes the measure worse}.
The inconsistent planner cannot justify adding the new route, even when the timetable is optimized.
We conclude that a consistent model can lead to good decisions even when the assumptions are wrong (shortest path versus logit), while an inconsistent model can lead to universally bad decisions.
It is remarkable that we already see this behavior in small examples, and we expect that the negative consequences of inconsistency may be worse in integrated models that allow bad decisions to propagate.

\begin{table}[t]
	\centering
	\begin{tabular}{lccc}
		\toprule
		~ & \multicolumn{3}{c}{Measure}\\
		\cmidrule(lr){2-4}
		\multicolumn{1}{l}{Timetable} & \multicolumn{1}{c}{Shortest Path} & \multicolumn{1}{c}{Logit $(\beta=0.1)$} & \multicolumn{1}{c}{Inconsistent} \\ \midrule
		Original (two routes)          & 30.00                       & 29.51                          & 31.42                                                                                                 \\
		Option 1: In-between        & 30.00                       & 28.23                          & 32.24                                                                                                 \\
		Option 2: Equidistant       & 31.67                       & 28.59                          & 33.40                                                                                                 \\
		Option 3: Optimal			& 29.44                       & 28.20                          & 31.87                                                                                                 \\ \hline
	\end{tabular}
	\caption{Three ways to evaluate the options in Figure~\ref{fig:scenarios}.}
	\label{tbl:measuresScen}
\end{table}
\section{Conclusions and Further Research}
\label{sec:con}

We presented a new framework for passenger routing and evaluating the service quality for route sets, timetables and line plans. We showed that the developed route choice models and measures are easy to interpret and can be computed efficiently, enabling practitioners to use the framework to better plan their networks.
The practical value of the new framework was demonstrated through several applications that were validated with real data from the Dutch and Swiss railway networks.
Most notably, it was found that the line plan measures are 99\% accurate in predicting the measure of the real system, effectively anticipating the impact of a timetable that is yet to be determined.
The new framework also generated several managerial insights about using the right model at the right time, when shortest path is all you need, and about the negative consequences of inconsistent models.
This includes the observation that, as one moves from immediate route choice to timetabling to line planning, the benefit of logit over shortest path diminishes, implying that the shortest path model, despite being slightly less accurate, may suffice for tactical or strategic planning purposes.

This paper lays the groundwork for a plethora of further research. The developed measures can be incorporated in virtually every public transport problem that concerns passengers such as line planning, timetabling, providing travel advice or determining tariff schemes. For problems with fixed route durations (e.g.\ line planning), the measures naturally lend themselves to decomposition schemes such as Benders decomposition due to their convex nature. Problems where the route durations depend on decision variables (e.g.\ timetabling) are more challenging, as this introduces non-convexity. This research also paves the way for new ways to evaluate and optimize the next generation of passenger transport systems, such as on-demand transport or shared mobility services. This requires non-trivial extensions to the framework, since these transport modes are not operated periodically. The application of more advanced choice models such as mixed logit or nested logit also seems highly relevant in this context, as these allow for individual or correlated preferences for modes and other heterogeneous factors. Finally, one could include aspects such as pricing, crowding or robustness. It would be interesting to study whether the insights found in this paper generalize to such settings. 

\footnotesize
\vspace{10pt}
\noindent\textbf{Acknowledgements} This research is partly supported by NSF LEAP-HI proposal NSF-1854684.
\normalsize

\DeclareRobustCommand{\VAN}[3]{#3}

\bibliography{references}

\appendix
\clearpage

\section{Proofs Section~\ref{sec:rs} (Route Sets)}
\label{app:properties_proofs}
This appendix provides the proofs for the properties in Table~\ref{tbl:measures}.
{
	\edef\thetable{\getrefnumber{tbl:measures}}
	\begin{table}[!h]
		\centering
		\measuretable{}
	\end{table}
	\addtocounter{table}{-1}
}

\paragraph*{Monotonicity}
With shortest path routing, increasing route lengths and removing routes can only increase $l_{\min}$, but is not required to do so strictly.
This implies (weak) monotonicity.
For logit routing with travel time evaluation, Figure~\ref{fig:measures_plots} provides a counterexample for monotonicity.
For logit routing with perceived travel time evaluation, strictly increasing a route length or removing a route strictly decreases the argument of the logarithm.
As $-\frac{1}{\beta} \log(.)$ is strictly decreasing, this strictly increases the measure as required.

\paragraph*{Consistency}
Shortest path routing minimizes travel time by definition, so this combination is consistent.
Logit routing assigns a positive probability to every route, and therefore does not minimize the travel time, which makes this combination inconsistent.
Next, we prove that logit routing minimizes perceived travel time, i.e., $\routingBold^{logit}$ minimizes $\mathcal{E}_{ptt}(\routes, .)$.
A similar proof is provided by \cite{anderson1988representative}.
We use the Lagrange multiplier rule to minimize $\mathcal{E}_{ptt}(\routes, p)$ over $p \in \mathcal{P}$.
The Lagrangian is given by
$$\mathcal{L}(\routes, p, \lambda) = \sum_{i \in \routes} l_i p_i + \frac{1}{\beta}\sum_{i\in \routes} p_i \log(p_i) + \lambda\left(1 - \sum_{i \in \routes} p_i\right), \textrm{ on domain } p \ge 0, \lambda \in \mathbb{R}.$$
Because the objective $\mathcal{E}_{ptt}(\routes, .)$ is convex and the constraint is linear, it is sufficient for optimality to find a probability vector $p^* \in P$ and a multiplier $\lambda^* \in \mathbb{R}$ such that the stationarity condition
$$\frac{d\mathcal{L}(\routes, p, \lambda)}{dp_i} = l_i +\frac{1}{\beta}\left(\log p_i +1\right) - \lambda = 0, \quad \forall i \in \routes$$
holds. It is straightforward to verify that this is the case for $$p^*_i = \frac{e^{-\beta l_i}}{\sum_{j\in \routes} e^{-\beta l_j}} = p_i^{logit} \textrm{ and } \lambda^* = \frac{1}{\beta}-\frac{1}{\beta}\log\left(\sum_{j\in \routes} e^{-\beta l_j}\right),$$
which proves consistency for logit routing with perceived travel time evaluation.
Finally, we prove that shortest path routing with perceived travel time evaluation is inconsistent through an example.
Consider route set $\routes=\{1,2\}$ with route lengths $l_1=1$, $l_2=2$ and parameter $\beta = 1$.
The shortest path routing is given by $p_1^{sp} = 1$, $p_2^{sp} = 0$, resulting in an evaluation of $\mathcal{E}(\routes, p^{sp}) = 1$.
However, the solution $p_1 = \frac{1}{2}$, $p_2=\frac{1}{2}$ provides a value of $\mathcal{E}(\routes, p) = \frac{3}{2} - \log(2) < 1$, which proves that this combination is inconsistent.

\section{Proofs Section~\ref{sec:tt} (Timetables)}
\label{app:tt_proofs}

\begin{proposition}[Constant Routing between Departures]
\label{prop:constantrouting}
	For shortest path routing and logit routing the route choice does not change between departures. That is, for route $i \in \routes$ and successor $\sigma(i) \in \routes$ it holds that $\routingBold(\routes_\timetable(t))$ is constant for $t \in (\timetable_i, \timetable_{\sigma(i)}]$.
\end{proposition}

\begin{proof}
	Until the next departure, all waiting times decrease by the same amount as time progresses.
	Subtracting a constant from all route lengths does not affect which length is the minimum, which proves the proposition for shortest path routing.
	For logit routing we have for any route $j\in \routes$ and constant $C \in \mathbb{R}$ that
	$$p_j^{logit} = \frac{e^{-\beta l_j}}{\sum_{k\in \routes} e^{-\beta l_k}} = \frac{e^{-\beta (l_j+C)}}{\sum_{k\in \routes} e^{-\beta (l_k+C)}},$$
	which completes the proof.
\end{proof}

\begin{proposition}[Translation Invariance]
\label{prop:trans_invar}
	Both $\mathcal{M}^{RS}_{sp}$ and $\mathcal{M}^{RS}_{logit}$ are translation invariant between departures.
	That is, for route $i \in \routes$ and successor $\sigma(i) \in \routes$ it holds that
	$$\mathcal{M}^{RS}(\routes_\timetable(t)) = [\theta_{\sigma(i)} - t]_T + \mathcal{M}^{RS}(\routes_\timetable(\timetable_{\sigma(i)})), \quad \forall t \in (\timetable_i, \timetable_{\sigma(i)}].$$
\end{proposition}

\begin{proof}
	By Proposition~\ref{prop:constantrouting}, the routing $p \in \mathcal{P}$ is constant on the interval $t \in (\timetable_i, \timetable_{\sigma(i)}]$.
	The route set measure on the interval is therefore given by the evaluation function for this fixed $p$.
	It follows that
	\begin{align*}
		\mathcal{M}^{RS}_{logit}(\routes_\timetable(t)) &= \mathcal{E}_{ptt}(\routes_\timetable(t), p)\\
		&= \sum_{j \in \routes} \left(l_j + [\timetable_j -t]_T\right) p_j + \frac{1}{\beta}\sum_{j\in \routes} p_j \log(p_j) &\textrm{(Definition $\mathcal{E}_{ptt}$ and $\routes_\timetable(t)$)}\\
		&= \sum_{j \in \routes} \left(l_j + [\timetable_j - \timetable_{\sigma(i)}]_T + [\timetable_{\sigma(i)} - t]_T \right) p_j + \frac{1}{\beta}\sum_{j\in \routes} p_j \log(p_j) &\textrm{(split the wait time)}\\
		&= [\timetable_{\sigma(i)} - t]_T + \sum_{j \in \routes} \left(l_j + [\timetable_j - \timetable_{\sigma(i)}]_T \right) p_j + \frac{1}{\beta}\sum_{j\in \routes} p_j \log(p_j)  &\textrm{($\sum_{j\in \routes} p_j = 1$)}\\
		&= [\timetable_{\sigma(i)} - t]_T + \mathcal{M}^{RS}_{logit}(\routes_\timetable(\timetable_{\sigma(i)})). &\textrm{(Definition $\mathcal{E}_{ptt}$ and $\routes_\timetable(t)$)}
	\end{align*}
	For $\mathcal{M}^{RS}_{sp}$ the evaluation function omits the term $\frac{1}{\beta}\sum_{j\in \routes} p_j \log(p_j)$, but the proof is identical.
\end{proof}

\begin{proposition}
	\label{prop:alg_char}
	Algorithm~\ref{alg:characteristic} correctly calculates the characteristic values $\delta, \tau$ and $\Delta$.
\end{proposition}
\begin{proof}
	The main steps of the algorithm are straightforward, as indicated in the main text.
	It remains to go over the edge cases and justify the use of the waiting time correction $w(i,j)$.
	The first edge case is when $\lvert \routes \rvert = 1$.
	With a single route, the regular definition of $\delta_1$ would evaluate to $[\timetable_i - \timetable{\pi(i)}]_T = 0$.
	Hence, the algorithm sets the correct value $\delta_1 = T$ and the corresponding $\tau_1 = l_1$ as a special case.
	The vertical jump evaluates to $\Delta_1 = \delta_1 + \tau_1 - \tau_1 = T$ as expected.
	
	Next, consider the waiting time correction
	\begin{equation}
		w(i,j) = \begin{cases}
			T & \textrm{ if } \timetable_i = \timetable_j,~j < i,\\
			[\timetable_j - \timetable_i]_T & \textrm{ else.}
		\end{cases}
	\end{equation}
	This correction ensures that route with the same departure time $\theta_i = \theta_j$ are processed sequentially according to their ordering, essentially as if there were small gaps between the departures.
	Recall that routes that depart at the same time are ordered according to their index (Section~\ref{sec:tt}).
	Consider processing route $i \in \routes$.
	If route $j \in \routes$ departs at the same time but is processed later $j \ge i$, then the regular waiting time $[\timetable_j - \timetable_i]_T = 0$ applies.
	If route $j$ has already been processed, then it will take time $T$ before $j$ is reached again, and hence a waiting time of $T$ is used.
	After all routes at this particular time have been processed, all waiting times have been increased by $T$.
	The correction also ensures that the total gap between departures $\sum_{i\in \routes} \delta_i$ adds up to $T$ as expected.
	This is straightforward for the regular case, but for the edge case where all routes $\routes=\{1,\hdots,n\}$ depart at the same time, the correction is necessary to ensure that the horizontal gap from $n$ back to $1$ is recorded as $T$.		
\end{proof}

\begin{corollary}
	\label{col:sum_T}
	The output of Algorithm~\ref{alg:characteristic} satisfies $\sum_{i\in \routes} \delta_i = T$ and $\sum_{i \in \routes} \Delta_i = T$.
\end{corollary}
\begin{proof}
	The first summation follows from the definition of $\delta_i$, and the discussion of the edge cases in the proof of Proposition~\ref{prop:alg_char}.
	The second summation follows from the identity $\Delta_i = \delta_{\sigma(i)} + \tau_{\sigma(i)} - \tau_i$ in Line~\ref{line:Delta} of Algorithm~\ref{alg:characteristic}, and the fact that the $\tau$-terms cancel in the summation.
	It follows that $\sum_{i \in \routes} \Delta_i = \sum_{i \in \routes} \delta_i = T$.
\end{proof}

\begin{proposition}
	\label{prop:tt_unif_general}
	Under uniform arrivals, the timetable measure for shortest path routing or logit routing can be calculated as
	\begin{equation}
		\mathcal{M}^{TT}(\routes, T, \timetableBold) = \frac{1}{T} \sum_{i\in \routes} \left( \frac{1}{2} \delta_i^2 + \tau_i \delta_i \right) = \frac{1}{T} \sum_{i\in \routes} \left( \frac{1}{2} \Delta_i^2 + \tau_i \Delta_i \right), \tag{\ref{eq:tt_uniform}}
	\end{equation}
	where the characteristic values $\delta$, $\tau$, $\Delta$ are obtained from Algorithm~\ref{alg:characteristic}.
\end{proposition}
\begin{proof}
	By definition of the expectation value, the timetable measure is obtained as the area under the observed route set measure function, divided by $T$.
	For the first expression, Figure~\ref{fig:delta_repres} provides a graphical proof: For a given route $i\in \routes$ the shaded area is the sum of the area of half a square ($\frac{1}{2} \delta_i^2$) and a rectangle ($\delta_i \tau_i$).
	If routes $i$, $\sigma(i)$ depart at the same time it follows that $\delta_{\sigma(i)} = 0$, which does not affect the argument.
	Summing over the areas and dividing by $T$ gives the result.
	The second expression is obtained by substituting $\Delta_i = \delta_{\sigma(i)} + \tau_{\sigma(i)} - \tau_i$ (Line~\ref{line:Delta}, Algorithm~\ref{alg:characteristic}) into the first expression:
	\begin{align*}
		\frac{1}{T} \sum_{i\in \routes} \left( \frac{1}{2} \delta_i^2 + \tau_i \delta_i \right) &= \frac{1}{T} \sum_{i\in \routes} \left( \frac{1}{2} \delta_{\sigma(i)}^2 + \tau_{\sigma(i)} \delta_{\sigma(i)} \right)\\
		&= \frac{1}{T} \sum_{i\in \routes} \left( \frac{1}{2} \left(\Delta_i - \tau_{\sigma(i)} + \tau_i\right)^2 + \tau_{\sigma(i)} \left(\Delta_i - \tau_{\sigma(i)} + \tau_i\right) \right)\\
		&= \frac{1}{T} \sum_{i\in \routes} \left( \frac{1}{2} \Delta_i^2 + \tau_i \Delta_i + \frac{1}{2}\left(\tau_i^2 - \tau_{\sigma(i)}^2\right)\right)\\
		&= \frac{1}{T} \sum_{i\in \routes} \left( \frac{1}{2} \Delta_i^2 + \tau_i\Delta_i \right).
	\end{align*}	
\end{proof}

\begin{proposition}
	\label{prop:tau_Delta_def}
	For logit routing, 	the output of Algorithm~\ref{alg:characteristic} satisfies
	\begin{equation*}
		\tau_i = l_i + \frac{1}{\beta} \log\left( \frac{1-e^{-\beta \Delta_i}}{1-e^{-\beta T}} \right) \quad  \forall i \in \routes.
	\end{equation*}
\end{proposition}
\begin{proof}
	By definition, the value $\Delta_i$ is the jump in the observed route set measure due to missing route $i \in \routes$.
	Before the jump, the value of the observed route set measure is $\tau_i$.
	Missing route $i \in \routes$ increases the waiting time for this route from $0$ to $T$.
	It follows from Line~\ref{line:tau_logit} that
	\begin{align*}
		\Delta_i &= -\frac{1}{\beta}\log \left(\sum_{j\in \routes} e^{-\beta \left(l_j + w(i,j)\right)} + e^{-\beta (l_i + T)} - e^{-\beta (l_i + 0)} \right) - \tau_i,\\
		&= -\frac{1}{\beta}\log \left(e^{-\beta \tau_i} + e^{-\beta l_i} \left(e^{-\beta T} - 1\right) \right) - \tau_i.
	\end{align*}	
	Note that $\Delta_i > 0$ by strict monotonicity of the route set measure.
	Take exponents on both sides and rearrange to obtain
	\begin{equation*}
		e^{-\beta \Delta_i} = \frac{e^{-\beta \tau_i} + e^{-\beta l_i} \left(e^{-\beta T} - 1\right)}{e^{-\beta \tau_i}} = 1 + e^{-\beta l_i} \frac{e^{-\beta T} - 1}{e^{-\beta \tau_i}},
	\end{equation*}
	\begin{equation*}
		e^{-\beta \tau_i} = e^{-\beta l_i} \frac{1 - e^{-\beta T}}{1 - e^{-\beta \Delta_i}}.
	\end{equation*}
	Finally, take logarithms to recover $\tau_i$:
	\begin{equation*}
		\tau_i = l_i - \frac{1}{\beta}\log\left(\frac{1 - e^{-\beta T}}{1 - e^{-\beta \Delta_i}}\right) = l_i + \frac{1}{\beta}\log\left(\frac{1 - e^{-\beta \Delta_i}}{1 - e^{-\beta T}}\right).
	\end{equation*}
\end{proof}

\begin{proposition}
	\label{prop:measure_logit_convex}
	The timetable measure for logit routing under uniform arrivals,
	\begin{equation}
		\mathcal{M}^{TT}_{logit}(\routes, T, \timetableBold) = \frac{1}{T} \sum_{i\in \routes} \left( \frac{1}{2} \Delta_i^2 + l_i \Delta_i + \frac{\Delta_i}{\beta}\log\left(\frac{1 - e^{-\beta \Delta_i}}{1 - e^{-\beta T}}\right) \right),\tag{\ref{eq:tt_logit_convex}}
	\end{equation}
	is strictly convex on the domain $\Delta \ge 0$ for parameters $\beta > 0$ and $T > 0$.
\end{proposition}
\begin{proof}
	Due to symmetry, it is sufficient to prove strict convexity for a single $\Delta_i$.
	To simplify the proof, we split the logarithm and remove the linear terms $l_i \Delta_i$ and $-\frac{\Delta_i}{\beta}\log\left(1 - e^{-\beta T}\right)$.
	Furthermore, we scale the variable with the transformation $x = \beta \Delta_i$ and scale the function by $\beta^2 T > 0$.
	None of these transformations affect strict convexity.
 	It remains to prove that
 	\begin{equation*}
 		f(x) = \beta^2T \frac{1}{T} \left(\frac{1}{2} \frac{x^2}{\beta^2} + \frac{x}{\beta^2} \log \left(1- e^{-x}\right)\right) = \frac{1}{2} x^2 + x \log \left(1-e^{-x}\right)
 	\end{equation*}
 	is strictly convex.
 	Interestingly, the term $x \log \left(1-e^{-x}\right)$ is not convex by itself.
 	
	We prove strict convexity by showing that the derivative
    \begin{align*}
 			f'(x) &= x + \frac{x}{1-e^{-x}} e^{-x} + \log\left(1-e^{-x}\right)\\
 			&= \frac{x}{1-e^{-x}} (1-e^{-x}) + \frac{x}{1-e^{-x}} e^{-x} + \log\left(1-e^{-x}\right) \\
 			&=  \frac{x}{1-e^{-x}} + \log\left(1-e^{-x}\right)
 	\end{align*}
 	is strictly increasing (or equivalently, the second derivative is positive on its domain), which is a sufficient condition for strict convexity.
	The first term is non-decreasing, as the derivative
	\begin{equation}
		\frac{d}{dx} \frac{x}{1-e^{-x}} = \frac{(1-e^{-x}) - x e^{-x}}{(1-e^{-x})^2} = \frac{e^{-x}}{(1-e^{-x})^2} (e^x - 1 - x) \ge 0.
	\end{equation}
	This follows from the well-known inequality $e^x \ge 1 + x$.
	The second term of $f'(x)$ is strictly increasing by inspection.
	It follows that $f'(x)$, the sum of these two terms, is strictly increasing, and therefore $f(x)$ is strictly convex.
\end{proof}

\begin{corollary}
	\label{col:infderiv}
	The derivative of the timetable measure for logit routing under uniform arrivals tends to $-\infty$ when $\Delta_i \rightarrow 0$ for any $i \in \routes$.
\end{corollary}
\begin{proof}
	The transformations applied in the proof above do not affect whether or not the limits at $\Delta_i \rightarrow 0$ are finite.
	For $x \rightarrow 0$, the first term of $f'(x)$ is finite, and the second term tends to $-\infty$.
	It follows that the derivative tends to $-\infty$ when $\Delta_i \rightarrow 0$ for any $i \in \routes$.
\end{proof}

\section{Proofs Section~\ref{sec:lp} (Line Plans)}
\label{app:lp_proofs}

\begin{proposition}
	\label{prop:match_sp}
	The output of the procedure \textsc{TimeTableShortestPath} (Algorithm~\ref{alg:construct}) is a feasible timetable for which the characteristic values $\delta$ and match the desired values set during the execution of the algorithm.
\end{proposition}
\begin{proof}
	The input satisfies $x_i\ge 0$ for every route $i \in \routes$.
	The timetable is constructed according to $\hat{\delta}_i \gets x_i$ and $\timetable_i \gets \timetable_{\pi(i)} + \hat{\delta}_i$.
	This timetable is feasible, as the gaps are non-negative and sum to $\sum_{i \in \routes} \hat{\delta}_i = \sum_{i \in \routes} x_i = T$ as expected.
	It is clear from the definition that when $\delta_i$ is calculated by Algorithm~\ref{alg:characteristic}, indeed $\delta_i = \hat{\delta}_i$ matches.
\end{proof}

\begin{lemma}
	\label{lemma:y_pos}
	Every optimal solution to Problem~\eqref{formulation:logit} satisfies $y > 0$.
\end{lemma}
\begin{proof}
	The statement is trivial for a single route, as $y_1 = T$ by constraint~\eqref{eq:logit:opt:total}.
	Now assume that there are at least two routes $\lvert \routes \rvert \ge 2$ and that $y_i = 0$.
	Due to \eqref{eq:logit:opt:total}, there must be a route $j \in \routes$, $j\neq i$, for which $y_j > 0$.
	By Corollary~\eqref{col:infderiv}, the derivative of the timetable measure tends to $-\infty$ as $y_i \rightarrow 0$.
	It follows that the solution can strictly be improved by adding a small amount $\varepsilon > 0$ to $y_i$ while subtracting it from $y_j$.
	This contradicts the assumption that $y$ is optimal, so it must be that $y_i > 0$.
\end{proof}

\begin{lemma}
	\label{lemma:logit_mu}
	The desired values set by \textsc{TimeTableLogit} (Algorithm~\ref{alg:construct}) satisfy $\hat{\tau}_i + \frac{\hat{\Delta}_i}{1-e^{-\beta \hat{\Delta}_i}} = \mu$ for every route $i \in \routes$, for some constant $\mu \in \mathbb{R}$.
\end{lemma}
\begin{proof}
	The procedure sets $\hat{\Delta}_i \gets y_i$ and $\hat{\tau}_i \gets l_i + \frac{1}{\beta}\log\left(\frac{1 - e^{-\beta \hat{\Delta}_i}}{1 - e^{-\beta T}}\right) = l_i + \frac{1}{\beta}\log\left(\frac{1 - e^{-\beta y_i}}{1 - e^{-\beta T}}\right)$, where $y$ is an optimal solution to Problem~\eqref{formulation:logit}.
	Lemma~\ref{lemma:y_pos} states that $y > 0$, which also implies $y < T$.
	Then it must be that the partial derivatives of objective~\eqref{eq:logit:opt:obj} are all equal: if the derivatives with respect to routes $i,j\in \routes$, $i \neq j$, are not equal, then a better solution is obtained by moving a small amount $\varepsilon$ between the routes without affecting feasibility.
	Hence, the value of
	\begin{align*}
		T \frac{d}{dy_i} \eqref{eq:logit:opt:obj} &= y_i + y_i \frac{1-e^{-\beta T}}{1-e^{-\beta y_i}} \frac{e^{-\beta y_i}}{1-e^{-\beta T}} + \left(l_i + \frac{1}{\beta}\log\left(\frac{1 - e^{-\beta \hat{\Delta}_i}}{1 - e^{-\beta T}}\right)\right)\\
		&= \hat{\Delta}_i + \hat{\Delta}_i \frac{e^{-\beta \hat{\Delta}_i}}{1-e^{-\beta \hat{\Delta}_i}} + \hat{\tau}_i = \hat{\tau}_i + \frac{\hat{\Delta}_i}{1-e^{-\beta \hat{\Delta}_i}}
	\end{align*}
	is the same for each route $i\in \routes$.
\end{proof}

\begin{lemma}
	\label{lemma:dontmiss}
	The desired values set by \textsc{TimeTableLogit} (Algorithm~\ref{alg:construct}) satisfy $\hat{\tau}_i + \hat{\Delta}_i > \hat{\tau}_j$ for every pair of routes $i, j \in \routes$.
\end{lemma}
\begin{proof}
	By Lemma~\ref{lemma:logit_mu}, it is equivalent to prove
	\begin{equation*}
		\mu- \frac{\hat{\Delta}_i}{1-e^{-\beta \hat{\Delta}_i}} + \hat{\Delta}_i > \mu - \frac{\hat{\Delta}_j}{1-e^{-\beta \hat{\Delta}_j}} \Longleftrightarrow \frac{\hat{\Delta}_i}{e^{\beta\hat{\Delta}_i} - 1} <  \frac{\hat{\Delta}_j}{1-e^{-\beta \hat{\Delta}_j}}.
	\end{equation*}
	Let $z_i = \beta \hat{\Delta}_i$ and $z_j = \beta \hat{\Delta}_j$.
	By Lemma~\ref{lemma:y_pos}, $\hat{\Delta} = y > 0$.
	Hence, it is sufficient to prove the following for $z_i > 0$, $z_j > 0$:
	\begin{equation*}
		\frac{z_i}{e^{z_i} - 1} <  \frac{z_j}{1-e^{-z_j}}.
	\end{equation*}
	Now use the well-known inequality $e^{z} > 1 + z$ for $z\neq 0$ as follows.
	On the left-hand side, the numerator and denominator are both positive.
	The inequality shows that the denominator is larger, and thus that the fraction is less than one.
	Similarly, the right-hand side has a larger numerator, such that the fraction is larger than one.
	This completes the proof.	
\end{proof}

\begin{proposition}
	\label{prop:match_logit}
	The output of the procedure \textsc{TimeTableLogit} (Algorithm~\ref{alg:construct}) is a feasible timetable for which the characteristic values $\Delta$ and $\tau$ match the desired values set during the execution of the algorithm.
\end{proposition}
\begin{proof}
	The procedure sets $\hat{\Delta}_i \gets y_i$, $\hat{\tau}_i \gets l_i + \frac{1}{\beta}\log\left(\frac{1 - e^{-\beta \hat{\Delta}_i}}{1 - e^{-\beta T}}\right) = l_i + \frac{1}{\beta}\log\left(\frac{1 - e^{-\beta y_i}}{1 - e^{-\beta T}}\right)$, and $\hat{\delta}_i \gets \hat{\tau}_{\pi(i)} + \hat{\Delta}_{\pi(i)} - \hat{\tau}_i$, where $y$ is an optimal solution to Problem~\eqref{formulation:logit}.
	Then the timetable is constructed according to $\timetable_i \gets \timetable_{\pi(i)} + \hat{\delta}_i$.
	It follows directly from Lemma~\ref{lemma:dontmiss} that $\hat{\delta} > 0$.
	Furthermore, we have $\sum_{i \in \routes} \hat{\delta}_i = \sum_{i \in \routes} \hat{\Delta}_i$ due to the telescoping sum.
	This sum equals $T$ as expected, due to constraint~\eqref{eq:logit:opt:total}.
	It follows that the timetable is feasible, and it is easy to see that indeed $\delta = \hat{\delta}$ matches when $\delta$ is calculated by Algorithm~\ref{alg:characteristic}.
	
	Next we prove that $\tau = \hat{\tau}$ matches when $\tau$ is calculated by Algorithm~\ref{alg:characteristic}.
	Without loss of generality, we prove the statement for $i=1$, and we relabel the other routes such that $i=1$, $\sigma(i) = 2$, etc.
	Note that all routes depart at different times ($\hat{\delta} > 0$), so no wait time correction is necessary.
	We obtain the wait time function $w(1,j) = \sum_{k=2}^j \delta_k = \hat{\tau}_1 + \sum_{k=1}^{j-1} y_k - \hat{\tau}_j$, which is due to the telescoping sum, $\hat{\delta} = \delta$ (as proven previously), and $\hat{\Delta} = y$ (by definition).
	
	Now calculate $\tau_1$ according to Algorithm~\ref{alg:characteristic}, with the goal to prove that it matches $\hat{\tau}_1$:
	\begin{align*}
		\tau_1 = -\frac{1}{\beta}\log\left(\sum_{j\in \routes} e^{-\beta (l_j + w(i,j))} \right)
		&= -\frac{1}{\beta}\log\left(\sum_{j\in \routes} e^{-\beta (l_j + \hat{\tau}_1 + \sum_{k=1}^{j-1} y_k - \hat{\tau}_j)} \right)\\
		&=\hat{\tau}_1 -\frac{1}{\beta}\log\left(\sum_{j\in \routes} e^{-\beta (l_j + \sum_{k=1}^{j-1} y_k - \hat{\tau}_j)} \right) \stackrel{?}{=} \hat{\tau}_1.
	\end{align*}
	As the $\hat{\tau}_i$s cancel, it remains to prove that the argument of the logarithm is equal to one:
	\begin{align*}
		\sum_{j\in \routes} e^{-\beta (l_j + \sum_{k=1}^{j-1} y_k - \hat{\tau}_j)} &= 1\\
		\sum_{j\in \routes} e^{-\beta \sum_{k=1}^{j-1} y_k} \left(\frac{1 - e^{-\beta y_j}}{1 - e^{-\beta T}}\right) &=1 &&\textrm{(by definition $\hat{\tau}_j$)}\\
		\sum_{j\in \routes} \left(e^{-\beta \sum_{k=1}^{j-1} y_k} - e^{-\beta \sum_{k=1}^{j} y_k}\right) &=1 - e^{-\beta T} &&\\
		e^{-\beta 0} - e^{-\beta \sum_{k=1}^{n} y_k} &=1 - e^{-\beta T}. &&\textrm{(telescoping sum)}
	\end{align*}
	This holds because $\sum_{k=1}^n y_k = T$ by constraint~\eqref{eq:logit:opt:total}.
	We conclude that indeed $\tau = \hat{\tau}$.
	
	It only remains to prove that $\Delta = \hat{\Delta}$ matches, given that $\delta = \hat{\delta}$ and $\tau = \hat{\tau}$ already match.
	Algorithm~\ref{alg:characteristic} constructs $\Delta_i$ according to $\Delta_i \gets \delta_{\sigma(i)} + \tau_{\sigma(i)} - \tau_i = \hat{\delta}_{\sigma(i)} + \hat{\tau}_{\sigma(i)} - \hat{\tau}_i$.
	By definition, Algorithm~\ref{alg:construct} constructs $\hat{\delta}_i \gets \hat{\tau}_{\pi(i)} + \hat{\Delta}_{\pi(i)} - \hat{\tau}_i$.
	This substitution cancels out $\hat{\tau}_{\sigma(i)}$ and $\hat{\tau}_i$ such that indeed $\Delta_i = \hat{\Delta}_i$.
	This completes the proof.	
\end{proof}

\begin{proposition}
	\label{prop:optimal_sp}
	The output of the procedure \textsc{TimeTableShortestPath} (Algorithm~\ref{alg:construct}) is an optimal timetable, and its measure matches the optimal objective value of Problem~\eqref{formulation:sp}.
\end{proposition}
\begin{proof}
	Consider an optimal timetable with measure $\mathcal{M}^*$ and preprocess it by removing all suboptimal routes.
	As explained in the main text, the resulting timetable satisfies $\tau_i = l_i$ for all routes $i\in \routes$ while the measure stays the same.
	A corresponding solution to Problem~\eqref{formulation:sp} is obtained by setting $x_i = 0$ for all routes that have been removed in preprocessing, and $x_i=\delta_i^*$ otherwise, where $\delta^*$ are the $\delta$-values of the optimal timetable.
	Feasibility follows directly from Corollary~\ref{col:sum_T}, and it is easy that this solution has an objective value~\eqref{eq:sp:opt:obj} that matches $\mathcal{M}^*$.
	We conclude that for the optimal objective value $v^*$, we have $v^* \le \mathcal{M}^*$.
	
	Now solve Problem~\eqref{formulation:sp} and use Algorithm~\ref{alg:construct} to construct a timetable based on this solution.
	Proposition~\ref{prop:match_sp} states that this is always possible, and that $\delta = x$.
	The resulting timetable is evaluated with~\eqref{eq:tt_uniform}, which assigns value $\frac{1}{T} \left(\frac{1}{2} \delta_i^2 + \tau_i \delta_i\right)$ $=\frac{1}{T} \left(\frac{1}{2} x_i^2 + \tau_i x_i\right)$ to each route $i \in \routes$.
	This value cannot be worse than the value assigned by Problem~\eqref{formulation:sp}, as the objective function~\eqref{eq:sp:opt:obj} instead uses $\frac{1}{T} \left(\frac{1}{2} x_i^2 + l_i x_i\right)$, with $\tau_i \le l_i$ by definition.
	We conclude that $\mathcal{M}$, the measure of the constructed timetable, satisfies $\mathcal{M} \le v^*$.
	Combining this with $\mathcal{M}^* \le \mathcal{M}$ (by definition) and $v^* \le \mathcal{M}^*$ (proven above) implies that $\mathcal{M} = \mathcal{M}^* = v^*$.
	That is, the constructed timetable is optimal, and its measure matches the optimal objective value.
\end{proof}

\begin{proposition}
	\label{prop:optimal_logit}
	The output of the procedure \textsc{TimeTableLogit} (Algorithm~\ref{alg:construct}) is an optimal timetable, and its measure matches the optimal objective value of Problem~\eqref{formulation:logit}.
\end{proposition}
\begin{proof}
	Given an optimal solution to Problem~\eqref{formulation:logit}, Proposition~\ref{prop:match_logit} states that a timetable can be constructed for which all the characteristic values match.
	As $\Delta = \hat{\Delta} = y$, the objective function~\eqref{eq:logit:opt:obj} is identical to the timetable measure~\eqref{eq:tt_logit_convex}.
	Hence, the timetable is optimal, and its measure matches the optimal objective value.
\end{proof}

\begin{proposition}
	\label{prop:xoptstruc}
	Every optimal solution to Problem~\eqref{formulation:sp} satisfies $x_i = \max\{0,\mu-l_i\}$ for every route $i\in \mathbb{R}$, for some constant $\mu \in \mathbb{R}$.
\end{proposition}
\begin{proof}
	Associate the multiplier $\mu/T$ with constraint~\eqref{eq:sp:opt:total} to create the Lagrangian:
	\begin{align*}
		\mathcal{L}(x, \mu/T) = \frac{1}{T} \sum_{i\in \routes} \left( \frac{1}{2} x_i^2 + x_i l_i \right) + \frac{\mu}{T} \left(T - \sum_{i \in \routes} x_i\right) = \frac{1}{T} \sum_{i\in \routes} \left( \frac{1}{2} x_i^2 + (l_i - \mu) x_i \right) + \mu.
	\end{align*}
	It follows from the KKT conditions that every optimal solution satisfies $x^*=\argmin_{x\ge 0} \mathcal{L}(x, \mu/T)$ for some $\mu \in \mathbb{R}$ \citep{Patriksson2008-SurveyContinuousNonlinear}.
	The Lagrangian is separable, and the optimum is found by independently minimizing $n$ convex parabolas over $\mathbb{R}_+$.
	It is elementary to show that the unique optimal solution is $x_i = \max\{0,\mu-l_i\}$.
\end{proof}

\begin{corollary}
	\label{col:xoptstruc:unique}
	Problem~\eqref{formulation:sp} has unique primal and dual optimal solutions.
\end{corollary}
\begin{proof}
	Constraint~\eqref{eq:sp:opt:total} requires $\sum_{i \in \routes} x_i = \sum_{i \in \routes} \max\{0,\mu-l_i\} = T$.
	It is elementary to show that this equation has a unique solution for $\mu$, which in turn defines a unique solution for $x$.
\end{proof}

\begin{corollary}
	\label{col:reverse}
	If the parameters in Problem~\eqref{formulation:sp} satisfy
	\begin{equation*}
		\sum_{i \in \routes} \left(L - l_i\right) = T,
	\end{equation*}
	for some route duration upper bound $L \ge \max_{j \in \routes}l_j$, then the unique optimal solution is given by $x_i=L - l_i$ for all routes $i \in \routes$.
\end{corollary}
\begin{proof}
	Choosing $\mu = L$ results in $x_i = \max\{0, L - l_i\} = L - l_i$, because $L$ is an upper bound.
	By assumption, $\sum_{i \in \routes} x_i = \sum_{i \in \routes} \left( L - l_i\right)= T$, making $\mu = L$ the unique dual optimal solution and $x_i = L-l_i$ the unique primal optimal solution.
\end{proof}

\begin{lemma}
	\label{lemma:output_delta_logit}
	The output of the procedure \textsc{TimeTableLogit} (Algorithm~\ref{alg:construct}) satisfies $$\delta_i = \frac{\Delta_i}{1-e^{-\beta \Delta_i}} - \frac{\Delta_{\pi(i)}}{1-e^{-\beta \Delta_{\pi(i)}}} e^{-\beta \Delta_{\pi(i)}}.$$
\end{lemma}
\begin{proof}
	By Proposition~\ref{prop:match_logit}, the output of the procedure is a feasible timetable, and all the characteristic values match the desired values.
	We obtain:
	\begin{align*}
		\delta_i &=
		\tau_{\pi(i)} + \Delta_{\pi(i)} - \tau_i &&\textrm{(Line~\ref{line:Delta}, Algorithm~\ref{alg:characteristic})}\\
		&= \mu - \frac{\Delta_{\pi(i)}}{1-e^{-\beta \Delta_{\pi(i)}}} + \Delta_{\pi(i)} - \mu + \frac{\Delta_i}{1-e^{-\beta \Delta_i}} &&\textrm{(Lemma~\ref{lemma:logit_mu})}\\
		&=- \frac{\Delta_{\pi(i)}}{1-e^{-\beta \Delta_{\pi(i)}}} e^{-\beta \Delta_{\pi(i)}} + \frac{\Delta_i}{1-e^{-\beta \Delta_i}}.
	\end{align*}
\end{proof}

\begin{lemma}
	\label{lemma:output_plogit}
	The output of the procedure \textsc{TimeTableLogit} (Algorithm~\ref{alg:construct}) satisfies $$p_n^{logit}(\routes_\timetable(\timetable_i)) = \frac{1-e^{-\beta \Delta_n}}{1-e^{-\beta T}} e^{-\beta \sum_{k=i}^{n-1} \Delta_k}.$$
\end{lemma}
\begin{proof}
	By Proposition~\ref{prop:match_logit}, the output of the procedure is a feasible timetable, and all the characteristic values match the desired values.
	Algorithm~\ref{alg:construct} sets the routes to depart in order of index, and route $n\in \routes$ is the last to depart before the next period starts.
	It follows that the waiting time between route $i \in \routes$ and route $n$ is given by $w(i,n) = \sum_{k=i+1}^n \delta_k$.
	Using the identity $\delta_i = \tau_{\pi(i)} + \Delta_{\pi(i)} - \tau_i$ (Line~\ref{line:Delta}, Algorithm~\ref{alg:characteristic}), the telescoping sum results in $w(i,n) = \tau_i + \sum_{k=i}^{n-1} \Delta_k - \tau_n$.
	We obtain:
	\begin{align*}
		p_n^{logit}(\routes_\timetable(\timetable_i))
		&= \frac{e^{-\beta (l_n + w(i,n))}}{\sum_{k\in \routes} e^{-\beta (l_k + w(i,k))}} && \textrm{(definition $p^{logit}$)}\\
		&= \frac{e^{-\beta (l_n + \tau_i + \sum_{k=i}^{n-1} \Delta_k - \tau_n)}}{e^{-\beta \tau_i}} && \textrm{(definition $w(i,n)$; definition $\tau_i$ (Line~\ref{line:tau_logit}, Alg.~\ref{alg:characteristic}))}\\
		&=e^{-\beta(l_n - \tau_n)} e^{-\beta \sum_{k=i}^{n-1} \Delta_k}\\
		&=\frac{1-e^{-\beta \Delta_n}}{1-e^{-\beta T}}e^{-\beta \sum_{k=i}^{n-1} \Delta_k}. && \textrm{(definition $\tau_n$ by Proposition~\ref{prop:tau_Delta_def})}
	\end{align*}
\end{proof}

\begin{proposition}
	\label{prop:face_prob}
	Uniform random travelers who face the optimal timetable created by\\\textsc{TimeTableLogit} (Algorithm~\ref{alg:construct}) choose route $i \in \routes$ with probability $\Delta_i/T$.
\end{proposition}
\begin{proof}
	Without loss of generality, we relabel the routes such that $i=n$, $\pi(i) = n-1$, etc., and we prove the statement for route $n$.
	The probability that a uniform random traveler arrives in the interval $t \in (\timetable_{\pi(i)}, \timetable_i]$ is equal to $\delta_i/T$ by definition.
	Within this interval, the routing is constant by Proposition~\ref{prop:constantrouting}.
	By the law of total probability, the probability that a random traveler selects route $n$ is then given by $\sum_{i \in \routes} \frac{\delta_i}{T} p_n^{logit}(\routes_\timetable(\theta_i))$.
	For convenience, we premultiply this quantity by a constant to obtain:
	\begin{align*}
		\mathrlap{\left(\frac{1-e^{-\beta T}}{1-e^{-\beta \Delta_n}} T\right)\sum_{i \in \routes} \frac{\delta_i}{T} p_n^{logit}(\routes_\timetable(\theta_i))} \quad\\
		&= \sum_{i\in \routes} \left(\frac{\Delta_i}{1-e^{-\beta \Delta_i}} - \frac{\Delta_{\pi(i)}}{1-e^{-\beta \Delta_{\pi(i)}}} e^{-\beta \Delta_{\pi(i)}}\right) e^{-\beta \sum_{k=i}^{n-1} \Delta_k} &&\textrm{(Lemmas~\ref{lemma:output_delta_logit} and \ref{lemma:output_plogit})} \\
		&\mathrlap{= \sum_{i\in \routes} \frac{\Delta_i}{1-e^{-\beta \Delta_i}} e^{-\beta \sum_{k=i}^{n-1} \Delta_k}- \sum_{i\in \routes} \frac{\Delta_{\pi(i)}}{1-e^{-\beta \Delta_{\pi(i)}}} e^{-\beta \left(\sum_{k=i}^{n-1} \Delta_k + \Delta_{\pi(i)}\right)} } \\
		&= \sum_{i\in \routes} \frac{\Delta_i}{1-e^{-\beta \Delta_i}} e^{-\beta \sum_{k=i}^{n-1} \Delta_k}- \sum_{i\in \routes} \frac{\Delta_i}{1-e^{-\beta \Delta_i}} e^{-\beta \left(\sum_{k=\sigma(i)}^{n-1} \Delta_k + \Delta_i\right)} &&\textrm{(shift summation index)} \\
		&= \sum_{i\in \routes} \frac{\Delta_i}{1-e^{-\beta \Delta_i}} \left(e^{-\beta \sum_{k=i}^{n-1} \Delta_k}- e^{-\beta \left(\sum_{k=\sigma(i)}^{n-1} \Delta_k + \Delta_i\right)}\right)\\
		&= \frac{\Delta_n}{1-e^{-\beta \Delta_n}} \left(e^{-\beta 0}- e^{-\beta \left(\sum_{k=1}^{n-1} \Delta_k + \Delta_n\right)}\right) &&\textrm{(difference 0 unless $i=n$)}\\
		&= \frac{\Delta_n}{1-e^{-\beta \Delta_n}} \left(1 - e^{-\beta T}\right) = \left(\frac{1-e^{-\beta T}}{1-e^{-\beta \Delta_n}} T\right) \frac{\Delta_n}{T}. &&\textrm{(Corollary~\ref{col:sum_T})}
	\end{align*}
	This proves that $\sum_{i \in \routes} \frac{\delta_i}{T} p_n^{logit}(\routes_\timetable(\theta_i)) = \Delta_n/T$.
	As route $n$ was chosen without loss of generality, this completes the proof for all routes.
\end{proof}

\begin{lemma}
	\label{lemma:active_spacing}
	Let the optimal solution $x$ to Problem~\eqref{formulation:sp} define an active route set $\mathcal{A} = \{i \in \routes \vert x_i > 0\}$.
	For a given timetable $\timetableBold$, let $\pi_{\mathcal{A}}(.)$ be the predecessor function restricted to route set $\mathcal{A}$.
	Every timetable $\timetableBold$ that is optimal under shortest path routing satisfies $\timetable_i - \timetable_{\pi_{\mathcal{A}}(i)} = x_i$ for all active routes $i \in \mathcal{A}$.	
\end{lemma}
\begin{proof}
	
	Consider an arbitrary optimal timetable $\timetableBold$.
	Following the process described in the proof of Proposition~\ref{prop:optimal_sp}, we use the timetable to construct an optimal solution $x$ to Problem~\eqref{formulation:sp}.
	That is, suboptimal routes are removed in preprocessing and are assigned the value $x_i = 0$.
	Let $\bar{\mathcal{A}}$ be the remaining set of routes and define $\pi_{\bar{\mathcal{A}}}(i)$ to be the predecessor function restricted to these routes.
	We continue following the process by setting $x_i = \timetable_i - \timetable_{\pi_{\bar{\mathcal{A}}}(i)}$ for al routes $i \in \bar{\mathcal{A}}$, and we conclude that $x$ is an optimal solution to Problem~\eqref{formulation:sp}.
	Note that $x_i > 0$ for all routes $i \in \bar{\mathcal{A}}$, as strictly suboptimal routes have been removed in preprocessing, and it is suboptimal for two identical routes to depart at the same time.
	By Corollary~\ref{col:xoptstruc:unique}, the optimal solution $x$ to Problem~\eqref{formulation:sp} is unique.
	Therefore, the constructed route set $\bar{\mathcal{A}}$ matches the active route set $\mathcal{A}$ obtained from solving Problem~\eqref{formulation:sp} directly.
	It follows that $\timetable_i - \timetable_{\pi_{\mathcal{A}}(i)} = x_i$ for all active routes $i \in \mathcal{A}$ for an arbitrary optimal timetable $\timetableBold$.
\end{proof}

\begin{definition}{\textsc{TimetablingMultipleDestinations}}\\
	Consider a periodic timetabling setting that serves travelers from a single origin to $m \in \mathbb{N}$ destinations within a period $T \in \mathbb{N}$. 
	Travelers are served by $n \in \mathbb{N}$ routes with route durations $l^j = (l_1^j, \hdots, l_n^j)$ that depend on their destination $j \in \{1, \hdots, m\}$.
	Given weights $w \in \mathbb{N}^m$ and a bound $B \in \mathbb{Q}$, does there exist a single timetable $\timetableBold \in \{0, 1, \hdots, T-1\}^n$ such that
	\begin{equation*}
		\sum_{j=1}^m w_j \mathcal{M}^{TT}_{sp}(l^j, T, \timetableBold) \le B?
	\end{equation*}
\end{definition}

\begin{proposition}
	\label{prop:nphard}
	\textsc{TimetablingMultipleDestinations} is NP-complete in the strong sense for fixed $m \ge 2$.
\end{proposition}
\begin{proof}
	
	\textsc{TimetablingMultipleDestinations} is in NP, as timetables can be evaluated in polynomial time with Algorithm~\ref{alg:characteristic} and Equation~\eqref{eq:tt_uniform}.
	We prove NP-completeness in the strong sense with a pseudo-polynomial reduction from the strongly NP-complete \textsc{3-Partition} problem \citep{garey1979computers}:  given a multiset \( S = \{s_1, s_2, \dots, s_{3k}\} \) of positive integers such that  $\Omega/4 < s_i < \Omega/2$ for all $i\in \{1, \hdots, 3k\}$ and $\sum_{i=1}^{3k} s_i = k \Omega$, can \( S \) can be partitioned into \( k \) disjoint triplets \( S_1, \dots, S_k \) such that each triplet sums to \( \Omega \)?
	
	Given an instance of \textsc{3-Partition}, we construct the following instance of \textsc{TimetablingMultipleDestinations} with a fixed number of $m=2$ destinations.
	The period is set to $T=2 k \Omega$ and we define $n=4k$ routes.
	The route durations are set to 
	\begin{align*}
		l^1 = (\underbrace{\Omega-s_1, \hdots, \Omega-s_{3k}}_{3k}, \underbrace{0, \hdots, 0}_k) &&\textrm{and}&& l^2 = (\underbrace{2\Omega, \hdots, 2\Omega}_{3k}, \underbrace{0, \hdots, 0}_k).
	\end{align*}
	We observe that, if Problem~\eqref{formulation:sp} were solved for each OD pair independently, we would obtain the unique optimal solutions $x^1 = (s_1, \hdots, s_{3k}, \Omega, \hdots, \Omega)$ and $x^2 = (0, \hdots, 0, 2\Omega, \hdots, 2\Omega)$.
	This follows from applying Corollary~\ref{col:reverse} with upper bounds $L=\Omega$ and $L=2\Omega$, respectively, as indeed $\sum_{i\in \routes} x^1_i = \sum_{i \in \routes} x^2_i = 2k\Omega = T$.
	Finally, the weights are set to $w_1=w_2=1$ and the bound $B$ is set to match the sum of the two independent minimum objective values, calculated by substituting $(l^1, T, x^1)$ and $(l^2, T, x^2)$ into Objective~\eqref{eq:sp:opt:obj}.
	As required for a pseudo-polynomial reduction, we observe that both the size and the values of the numerical parameters of the new instance are polynomially bounded by the size and numerical values of the original instance.
	The bound $B \in \mathbb{Q}$ can be encoded as a fraction of two polynomially-bounded values, which follows immediately from the formula of Objective~\eqref{eq:sp:opt:obj}.
	
	We conclude the proof by showing that the constructed \textsc{TimetablingMultipleDestinations} instance is a \texttt{YES}-instance if and only if the original \textsc{3-Partition} instance is a \texttt{YES}-instance:
	\begin{itemize}
		\item $\Longrightarrow$: If the original \textsc{3-Partition} instance is a \texttt{YES}-instance, then there exist disjoint triplets $S_1, \hdots, S_k$ such that each triplet sums to $\Omega$.
		Next, we reorder the routes as follows: all routes in $S_1$, route $3k+1$, all routes in $S_2$, route $3k+2$, $\hdots$, all routes in $S_k$, route $4k$.
		Finally, we use Algorithm~\ref{alg:construct} with input $x^1$ to construct a timetable that is optimal for the first OD pair.
		By construction, the routes $3k+1$, $\hdots$, $4k$ are evenly spaced $2\Omega$ apart.
		This follows from the fact that the $x^1$-values of each triplet progress the timetable by a total of $\Omega$, and the $x^1$ value of routes $i \in \{3k+1, \hdots, 4k\}$ progresses the timetable by another $\Omega$.
		For the second OD pair, after removing the unused routes $i \in \routes$ for which $x^2_i = 0$, this matches the optimal spacing defined by $x^2$.
		It follows that this timetable allows both OD pairs to attain their independent minimum objective value, which implies that \textsc{TimetablingMultipleDestinations} instance is a \texttt{YES}-instance.

		\item $\Longleftarrow$: If the constructed \textsc{TimetablingMultipleDestinations} instance is a \texttt{YES}-instance, then there exists a timetable $\timetableBold$ that allows both OD pairs to attain their independent minimum objective value.
		Applying Lemma~\ref{lemma:active_spacing} to optimal solution $x^2$ shows that $\timetableBold$ must have routes $3k+1$, $\hdots$, $4k$ at equidistant intervals $2\Omega$ apart.
		Applying the same lemma to optimal solution $x^1$ dictates that the gap before each route $i \in \routes$ must be equal to $x^1_i$.
		For the equidistant routes $i \in \{3k+1, \hdots, 4k\}$, this value is $x^1_i = \Omega$, leaving $k$ gaps of size $\Omega$ for the other routes.
		The remaining routes $i \in \{1, \hdots, 3k\}$ have values $x^1_i = s_i$, with $\Omega/4 < s_i < \Omega/2$ by definition.
		It follows that each gap of size $\Omega$ can only be filled by a triplet of values from $S$ that sum to $\Omega$.
		Thus, a valid 3-partition of $S$ must exist, which proves that the \textsc{3-Partition} instance is a \texttt{YES}-instance.
	\end{itemize}
\end{proof}

\end{document}